\documentclass{amsart}

\usepackage{graphicx}
\usepackage{natbib}
\usepackage{amssymb}
\usepackage{amsthm}
\usepackage{amsmath}
\usepackage{mathtools}
\newtheorem{theorem}{Theorem}

\newtheorem{lemma}{Lemma}

\newtheorem{definition}{Definition}
 
\usepackage{multicol}
\usepackage{lineno}
\usepackage{array}
\usepackage{boldline}
\usepackage{hyperref}
\usepackage[nosort]{cleveref}
\usepackage{tikz}
\usetikzlibrary{shapes.geometric, arrows}
\usetikzlibrary{er,positioning}
\usetikzlibrary{shapes}
\usepackage{subcaption}
\usepackage[font=small,skip=5pt]{caption}

\newcommand*\circled[1]{\tikz[baseline=(char.base)]{%
		\node[shape=circle,fill=white!80!black,draw,minimum size=20] (char) {#1};}}

\usepackage{algorithm}
\usepackage{algpseudocode}

\algnewcommand\True{\textbf{True}\space}
\algnewcommand\False{\textbf{False}\space}
\algloopdefx{NFor}[1]{\textbf{for all} #1 \textbf{do}}
\algloopdefx{NIf}[1]{\textbf{if} #1 \textbf{then}}

\usepackage{multirow}
\usepackage{lscape}
\usepackage{adjustbox}

\usepackage{url}
\usepackage{doi}

\DeclareMathOperator{\dom}{dom}
\DeclareMathOperator{\opt}{opt}
\DeclareMathOperator*{\argmax}{arg\,max}

\newcommand{\linprogref}[1]{\textup{(#1)}}

\title[Improving and benchmarking of algorithms for decision making]{Improving and benchmarking of algorithms for decision making with lower previsions}

\author{Nawapon Nakharutai}
\address{Durham University, Department of Mathematical Sciences, UK}
\email{nawapon.nakharutai@gmail.com}

\author{Matthias C. M. Troffaes}
\address{Durham University, Department of Mathematical Sciences, UK}
\email{matthias.troffaes@durham.ac.uk}

\author{Camila C. S. Caiado}
\address{Durham University, Department of Mathematical Sciences, UK}
\email{c.c.d.s.caiado@durham.ac.uk}


\begin{document}

\begin{abstract}
Maximality, interval dominance, and E-admissibility are three well-known criteria for decision making under severe uncertainty using lower previsions. We present a new fast algorithm for finding maximal gambles. We compare its performance to existing algorithms, one proposed by Troffaes and Hable (2014), and one by Jansen, Augustin, and Schollmeyer (2017). To do so, we develop a new method for generating random decision problems with pre-specified ratios of maximal and interval dominant gambles.

Based on earlier work, we present efficient ways to find common feasible starting points in these algorithms. We then exploit these feasible starting points to develop early stopping criteria for the primal-dual interior point method, further improving efficiency. We find that the primal-dual interior point method works best.

We also investigate the use of interval dominance to eliminate non-maximal gambles. This can make the problem smaller, and we observe that this benefits Jansen et al.'s algorithm, but perhaps surprisingly, not the other two algorithms. We find that our algorithm, without using interval dominance, outperforms all other algorithms in all scenarios in our benchmarking.
\end{abstract}

\maketitle

\section{Introduction}\label{intro}

Consider a subject who needs to choose from a set of possible decisions. Each decision leads to an uncertain reward, depending on her decision and on the state of nature revealed after the decision.
The reward could be anything, for example, money, food, or a lottery ticket.
For simplicity, we will assume that rewards are expressed on a utility scale.
In this way, we can view an uncertain reward, and thereby, each decision, as a bounded real-valued function on the set of states of nature. Such function will be called a \emph{gamble}.

The subject wants to choose gambles that lead to the best possible reward. We assume that, for any set of gambles, the subject can identify a subset of gambles that she does not want to choose. We say that a gamble is \emph{optimal} in a given set of gambles if the subject is not committed to eliminate it.

If we assume that, under the usual rationality assumptions,  the subject can specify a complete probability measure on the state of nature, then she should simply choose a gamble that maximises her expected utility \citep{anscombe:1963}. However, when only little information is available, the subject may not be able to specify a complete probability measure. In that case, the subject might consider other ways to express her beliefs. In this study, we will assume that the subject is able to model her beliefs using a \emph{lower prevision}, or equivalently, through probability bounding \citep{1975:williams:condprev,2007:williams:condprev}.
\emph{Maximality} \citep[\S 3.9.1-3.9.3, pp. 160-162]{1991:walley} and \emph{interval dominance} \citep[\S 2.3.3, pp.68-69]{2003:zaffalon::dementia} are well-known decision criteria induced by strict partial orders associated with lower previsions \citep{2007:troffaes:decision:intro}.

Several authors proposed algorithms for finding maximal gambles, for example, \citet{2011:Kikuti:Cozman:Filho}, \citet[p.~336]{2014:troffaes:itip:computation} and \citet{2017:Jansen:Augustin:Schollmeyer}.

\Citet[Algorithm~16.4, p.~336]{2014:troffaes:itip:computation} present an incremental algorithm where once some maximal gambles in the sets are known, we should compare the remaining gambles against those maximal gambles first. Additionally, \citet[p.~336]{2014:troffaes:itip:computation} suggest that sorting all gambles in advance, e.g. by expectation, could help the algorithm to perform better. In this paper, we propose a new algorithm that incorporates this suggestion, and we confirm that this leads to a considerable speed-up.

In the algorithms proposed in \citet{2011:Kikuti:Cozman:Filho} and \citet[p.~336]{2014:troffaes:itip:computation}, to check maximality of each gamble, one has to evaluate the sign of the lower prevision of several gambles.
This can be done by solving several linear programming problems \citep[p.~331]{2014:troffaes:itip:computation}. From earlier work, we know that the primal-dual interior point method is particularly suitable when working with lower previsions \citep{2018:Nakharutai:Troffaes:Caiado}. We propose early stopping criteria to determine more quickly the sign of lower previsions, exploiting the fact that primal-dual methods solve both the primal and dual simultaneously. We also use results from \citet{2018:Nakharutai:Troffaes:Caiado} to quickly obtain feasible starting points, further improving the efficiency of these methods.

\Citet{2017:Jansen:Augustin:Schollmeyer} proposed an algorithm that can verify whether a gamble is maximal by solving a single larger linear program. To improve this algorithm, we exploit the fact that if a gamble is not maximal in a given iteration, then it can be excluded from all future iterations. We verify that this improved version performs slightly faster.

As all maximal gambles are also interval dominant, \citet{2007:troffaes:decision:intro} suggested to eliminate non-maximal gambles by applying interval dominance first. When most gambles are not interval dominant, this can eliminate many non-maximal gambles early on. This might be useful, because interval dominance is easier to check. In this paper, we will compare the above mentioned algorithms for finding maximal gambles, with and without applying interval dominance.

The contributions of the paper are as follows.
We propose a new algorithm for finding maximal gambles and compare its performance to the two algorithms proposed by \citet[p.~336]{2014:troffaes:itip:computation} and \citet{2017:Jansen:Augustin:Schollmeyer}. For benchmarking, we propose an algorithm for generating sets of gambles which have a precisely given number of maximal and interval dominant gambles. 
For the algorithm in \citet[p.~336]{2014:troffaes:itip:computation}, and also for our new algorithm, we solve a sequence of linear programs by the primal-dual method, because we can easily find a common feasible starting point. For further improvement, we also develop early stopping criteria, so the method can stop iterating before achieving an optimal solution.

The paper is organised as follows. In \cref{sec:decision}, we present the basic concepts used in the remained of the paper. First, we give a brief outline of lower previsions and natural extension. Then, we present several decision criteria for lower previsions.  In \cref{sec:algorithm}, we discuss several algorithms from the literature, and propose a new algorithm for finding maximal gambles.  To benchmark these algorithms, in \cref{sec:benchmark}, we provide an algorithm for generating sets of gambles, and compare the performance of various decision algorithms on generated sets. Proofs of technical results from this section are provided in \ref{pf:appendix}. \Cref{sec:conclusion} concludes the paper.

\section{Decision making with lower previsions}\label{sec:decision}

In this section, we first explain lower previsions and natural extension. Then, we present three decision criteria: maximality, E-admissibility, and interval dominance.
For more about the relation between these criteria and their advantages and disadvantages, we refer to \citet{2014:troffaes:itip:decision} and \citet{2007:troffaes:decision:intro}.

\subsection{Lower previsions}

Let $\Omega$ denote the set of states of nature, and let $\mathcal{L}$ denote the set of all gambles (i.e. bounded real-valued functions) on $\Omega$.
A \emph{lower prevision} $\underline{P}$ is a real-valued function defined on some subset of $\mathcal{L}$. We denote the domain of $\underline{P}$ by $\dom\underline{P}$. Given a gamble $f\in\dom\underline{P}$, we interpret $\underline{P}(f)$ as the subject's supremum buying price for $f$,
i.e. for all $\alpha<\underline{P}(f)$,
she is willing to accept the gamble $f-\alpha$.
It has been argued that this is a good way for the subject to model
her uncertainty about $\Omega$, especially under severe uncertainty
\citep{1991:walley,2008:miranda::survey:lowprevs,2014:miranda:itip:lowerprevision,2014:troffaes:decooman::lower:previsions}.

We say $\underline{P}$ \emph{avoids sure loss} if 
for all $ n \in \mathbb{N}$, all  $\lambda_{1}, \dots,\lambda_{n} \geq 0$, and all $f_{1}, \dots,f_{n} \in \dom\underline{P}$, it holds that \citep[p.~42]{2014:troffaes:decooman::lower:previsions}:
\begin{equation}\label{eq:3.1}
\max_{\omega\in \Omega} \left( \sum_{i=1}^{n} \lambda_{i}\left[f_{i}(\omega)-\underline{P}(f_{i})\right] \right) \geq 0. 
\end{equation} 
If $\underline{P}$ does not avoid sure loss, then there is a finite combination of gambles $f_{1}, \dots,f_{n} \in \dom\underline{P}$ such that for some $\lambda_{1}, \dots,\lambda_{n} \geq 0$: 
\begin{equation}
\max_{\omega\in \Omega} \left(\sum_{i=1}^{n} \lambda_{i}f_{i}(\omega) \right)  < \sum_{i=1}^{n} \lambda_{i}\underline{P}(f_{i})
\end{equation}
which means that the subject is willing to pay more than what she could ever gain, which does not make sense \citep[p.~44]{2014:troffaes:decooman::lower:previsions}. Therefore, throughout this study, we assume that all lower previsions avoid sure loss.
 
The conjugate upper prevision $\overline{P}$ on $-\dom \underline{P}\coloneqq\{-f: f\in \dom \underline{P}\}$ is defined by $\overline{P}(f)\coloneqq -\underline{P}(-f)$. It represents the subject's infimum selling price for $f$ \citep[p.~41]{2014:troffaes:decooman::lower:previsions}. 

We can extend $\underline{P}$ to the set of all gambles $\mathcal{L}$ via its \emph{natural extension} $\underline{E}$. For any gamble $g$, the natural extension $\underline{E}(g)$ corresponds to the supremum price a subject should be willing to pay for $g$, given the prices $\underline{P}(f)$ for all $f \in \dom \underline{P}$ \citep[p.~47]{2014:troffaes:decooman::lower:previsions}.

\begin{definition} \citep[p.~47]{2014:troffaes:decooman::lower:previsions}
Let $\underline{P}$ be a lower prevision. The natural extension of $\underline{P}$ is defined on all $f \in \mathcal{L}$ by:
\begin{equation}
\underline{E}(f)\coloneqq\sup \left\lbrace\alpha \in \mathbb{R}\colon f -\alpha \geq \sum_{i=1}^{n} \lambda_{i}(f_{i}-\underline{P}(f_{i})) , n\in\mathbb{N},\,f_{i} \in  \dom\underline{P},\,\lambda_{i} \geq 0\right\rbrace.
\end{equation} 
\end{definition}

Note that $\underline{E}$ is finite, and therefore, is a lower prevision, if an only if $\underline{P}$ avoids sure loss \citep[p.~68]{2014:troffaes:decooman::lower:previsions}. In the case that both $\Omega$ and $\dom \underline{P}$ are finite, for any gamble $f$, $\underline{E}(f)$  can be calculated by solving a linear program \citep[p.~331]{2014:troffaes:itip:computation}. 

We denote the conjugate of $\underline{E}$ by $\overline{E}$ which is given by
\begin{align}
  \overline{E}(f)& \coloneqq -\underline{E}(-f)\\
 & = \inf \left\lbrace\beta \in \mathbb{R}\colon \beta - f \geq \sum_{i=1}^{n} \lambda_{i}(f_{i}-\underline{P}(f_{i})), n \in \mathbb{N}, f_{i} \in \dom\underline{P}, \lambda_{i} \geq 0\right\rbrace.
\end{align}

\subsection{Decision criteria}

We first define two strict partial orderings on $\mathcal{L}$, and then define optimality through maximality with respect to either of these two strict partial orderings.

\begin{definition}
For any two gambles $f$ and $g$, we say that $f \succ g$ whenever 
\begin{equation}
\underline{E}(f-g) > 0.
\end{equation}
\end{definition}

Note that \citet[\S 3.8.1]{1991:walley} uses a stronger ordering, which also includes pointwise dominance. Here, we follow 
\citet{2007:troffaes:decision:intro},
\citet[\S 16.3.2]{2014:troffaes:itip:computation} and \citet{2017:Jansen:Augustin:Schollmeyer}  and simply omit pointwise dominance from our definition.

\begin{definition}\citep[p.~194]{2014:troffaes:itip:decision}
For any two gambles $f$ and $g$, we say that $f\sqsupset g$ whenever 
\begin{equation}
\underline{E}(f) > \overline{E}(g).
\end{equation}
\end{definition}

Given any strict partial order on $\mathcal{L}$,
we can define a notion of optimality
through maximality with respect to that order:
\begin{definition}
Let $\pmb{>}$ be a strict partial order on $\mathcal{L}$, and let $\mathcal{K}$ be a finite subset of $\mathcal{L}$. The \emph{set of maximal gambles in $\mathcal{K}$ with respect to $\pmb{>}$} is then defined by:
\begin{equation}\label{eq:opt_max}
\opt_{\pmb{>}}(\mathcal{K}) \coloneqq \{f \in \mathcal{K}\colon (\forall g \in \mathcal{K})(g \not \pmb{>} f )\}.
\end{equation}
\end{definition}

We call $\opt_{\succ}(\mathcal{K})$ the set of \emph{maximal} gambles in $\mathcal{K}$ and $\opt_{\sqsupset}(\mathcal{K})$ the set of \emph{interval dominant} gambles in $\mathcal{K}$.

Finally, we also need to define E-admissibility, which is yet
another decision criterion.
First, we need some notation.
The unit simplex is the set of all probability mass functions:
\begin{equation}
\Delta \coloneqq \left\{ p \in \mathbb{R}^\Omega\colon p \geq 0\text{\ and }\sum_{\omega \in \Omega} p(\omega) = 1 \right\}.
\end{equation}
The \emph{credal set} of a lower prevision $\underline{P}$ is defined by \citep[p.37]{2014:miranda:itip:lowerprevision}:
\begin{equation}
\mathcal{M} \coloneqq \{p \in \Delta\colon \forall f \in \dom\underline{P},\ E_p(f) \geq\underline{P}(f) \}.
\end{equation}
\begin{definition}[E-admissibility]\citep[p.~336]{2014:troffaes:itip:computation} 
A gamble $f$ is \emph{E-admissible} in $\mathcal{K}$  if there is $p \in \mathcal{M}$ such that 
\begin{align}
\forall g\in \mathcal{K} \colon E_p(f) \geq E_p(g).
\end{align}
\end{definition}
The set of all E-admissible gambles in $\mathcal{K}$
is denoted by $\opt_{\mathcal{M}}(\mathcal{K})$.

Note that \citep{2007:troffaes:decision:intro}:
\begin{equation}
\opt_{\mathcal{M}}(\mathcal{K})\subseteq\opt_{\succ}(\mathcal{K})  \subseteq \opt_{\sqsupset}(\mathcal{K}).
\end{equation}
If a gamble is not interval dominant, then it is not maximal.
Consequently, if there are many gambles in the set, one may want to eliminate non-maximal gambles in $\mathcal{K}$ by applying interval dominance first \citep{2007:troffaes:decision:intro}.

Similarly, if we find an E-admissible gamble $f$ in $\mathcal{K}$, then $f$ is immediately maximal   \citep[\S 3.9.4]{1991:walley}.
In \cref{sec:algorithm}, we will show how we can quickly find an E-admissible gamble in $\mathcal{K}$ to speed up algorithms for finding $\opt_{\succ}(\mathcal{K})$.

\section{Improving algorithms for finding maximal gambles}\label{sec:algorithm}

\subsection{Algorithms for finding maximal gambles}

In this section, we will discuss algorithms for finding $\opt_{\succ}(\mathcal{K})$. We study two algorithms from the literature, and propose a new algorithm based on a suggestion from \citet[p.~336]{2014:troffaes:itip:computation}.

One can see that a gamble $f$ is maximal in $\mathcal{K}$ only if  
\begin{equation}
\forall g \in \mathcal{K}: \overline{E}(f-g) \geq 0.
\end{equation}
Suppose that there are $m$ possible outcomes in $\Omega$, $k$ gambles in $\mathcal{K}$ and $n$ gambles in $\dom\underline{P}$ where $\underline{P}$ avoids sure loss.
To check whether $f$ is maximal in $\mathcal{K}$, we have to calculate $\overline{E}(f-g)$ for all $g \in \mathcal{K}\setminus\{f\}$. Let $h\coloneqq g-f$, we can calculate $\underline{E}(h)$ through either \linprogref{P1} or \linprogref{D1}:
\begin{align}
\label{P1:1}\tag{P1a}
\linprogref{P1} &&
\min \quad & \sum_{\omega \in \Omega} h(\omega)p(\omega)  \\
\label{P1:2}\tag{P1b}
&& \text{subject to}\quad & \forall g_i\in\dom\underline{P} \colon \sum_{\omega \in \Omega} (g_i(\omega)-\underline{P}(g_i))p(\omega) \geq 0\\
\label{P1:3}\tag{P1c}
&& & \sum_{\omega \in \Omega}p(\omega) = 1\\
\label{P1:4}\tag{P1d}
&&\text{where} \quad  & \forall \omega\colon p(\omega) \geq 0,
\end{align}

\begin{align}
\tag{D1a}
\linprogref{D1}  &&
\max\quad & \alpha \\
\label{thm1:5}\tag{D1b}
&& \text{subject to}\quad & \forall \omega \in\Omega\colon  \sum_{i=1}^{n} (g_i(\omega)-\underline{P}(g_i))\lambda_{i} + \alpha \leq h(\omega) \\
\label{thm1:3}\tag{D1c}
&& \text{where} \quad & \forall i\colon \lambda_{i} \geq 0 \quad (\alpha \text{ free}).
\end{align}
$\underline{E}(h)$ is precisely the optimal value of \linprogref{P1} (or \linprogref{D1}). The problem \linprogref{D1} is an unconditional case of the linear program in \cite[p.~331]{2014:itip}. Note that  \linprogref{P1} has $n+1$ constraints and $m$ variables. So, to determine all maximal gambles, we must solve (at most) $k(k-1)$ of these linear programs.

\Citet[algorithm 16.4, p.~336]{2014:troffaes:itip:computation} proposed the following strategy for finding maximal gambles: once a non-maximal gamble is detected, it is no longer compared with other gambles. Indeed, if $f$ is non-maximal, then there will be some gamble $g$ that dominates $f$. However, if $g$ dominates $f$, and $f$ dominates $h$, then 
$g$ will also dominates $h$ as well. Consequently, every non-maximal gamble is dominated by at least one maximal gamble in $\mathcal{K}$. Therefore, the algorithm no longer needs to consider non-maximal gambles as soon as they are deemed non-maximal (see \cref{alg:FindMax:matt}).
\begin{algorithm}
\caption{Find the set of maximal gambles in $\mathcal{K}$}\label{alg:FindMax:matt}
\begin{algorithmic}[1]
\Require a set of $k$ gambles $\mathcal{K} = \{ f_1,\dots, f_k\}$
\Ensure an index set of $\opt_{\succ}(\mathcal{K})$
\State $I \gets \emptyset$  \Comment{an index set of $\opt_{\succ}(\mathcal{K})$}
\For{$i=1:k$}
\If{IsNotDominated1$(I,i)$} 
\State{\indent $I\gets I \cup \{i\}$}\Comment{$f_i$ is maximal}
\EndIf
\EndFor
\State \Return $I$
\State{\textbf{where} IsNotDominated1($I,i$)}
\State{\indent \textbf{for} $j\in I \cup \{i+1,\dots,k\}$} 
\State{\indent \indent \textbf{if} $\underline{E}(f_j-f_i) > 0$  \textbf{then return False}} \Comment{$f_i$ is dominated by $f_j$}
\State{\indent\Return \True}
\end{algorithmic}
\end{algorithm}

For \cref{alg:FindMax:matt}, if the first considered gamble happens
to be the only maximal gamble in $\mathcal{K}$, then the algorithm
only needs to solve $2(k-1)$ linear programs
\citep[p.~336]{2014:troffaes:itip:computation}. Specifically, to
verify that none of gambles in the set dominate the first gamble, the
algorithm first needs to solve $k-1$ linear programs. Next, for each of
the remaining gambles, the algorithm compares it with the only existing
maximal gamble, so the algorithm additionally solves $k-1$ linear
programs.  If all gambles in $\mathcal{K}$ are maximal, then the
method needs to solve $k(k-1)$ linear programs, because for each
gamble, the algorithm has to compare it with the other $k-1$ gambles
\citep[p.~335]{2014:troffaes:itip:computation}.

This shows that we could speed up \cref{alg:FindMax:matt} by early identifying some maximal gambles in $\mathcal{K}$, for example, via E-admissibility. Specifically, if a gamble $f$ is E-admissible in $\mathcal{K}$, then $f$ is also maximal in $\mathcal{K}$ \citep[p.~196]{2014:troffaes:itip:decision}. We can simply find one of E-admissible gambles as follows.
We first sort all gambles in $\mathcal{K}$ as $f_1,\dots,f_k$ such that for some $p \in \mathcal{M}$, for all $j >i$:
\begin{equation}\label{eq:E_p(f_j-f_i)}
E_p(f_j) - E_p(f_i) \geq 0.
\end{equation}
Then $f_k$ has the highest expectation, and therefore $f_k$ is E-admissible in $\mathcal{K}$. We can then improve \cref{alg:FindMax:matt} by initially setting $\opt_{\succ}(\mathcal{K}) = \{f_k \}$. 
To obtain such $p \in \mathcal{M}$, we can solve \linprogref{P1} with $h=0$.
In principle, one could identify further maximal gambles by finding all E-admissible gambles, for example by using one of the algorithms proposed by \citet{Kikuti05partiallyordered} or \citet{Utkin2005PowerfulAF}.
However, these algorithms require solving $k$ linear programs where $k=|\mathcal{K}|$,
and they are more complex than \linprogref{P1} with $h=0$.
Therefore, we do not apply those algorithms here.

In addition to identifying one E-admissible gamble in $\mathcal{K}$, sorting gambles with respect to the expectation also saves a lot of comparison steps in \cref{alg:FindMax:matt} for finding $\opt_{\succ}(\mathcal{K})$. Specifically, to determine whether gamble $f_i$ is maximal in $\mathcal{K}$, we need to evaluate only $\underline{E}(f_j-f_i)$ when $j >i$, because we immediately know that
\begin{equation}
 \forall i\le j\colon \underline{E}(f_i-f_j)\leq E_p(f_i-f_j) \leq 0.
\end{equation}
An algorithm for finding maximal gambles that exploits sorting gambles is presented in \cref{alg:FindMax2}.

\begin{algorithm}
\caption{Find the set of maximal gambles in $\mathcal{K}$}\label{alg:FindMax2}
\begin{algorithmic}[1]
\Require a set of $k$ gambles $\mathcal{K} = \{ f_1,\dots, f_k\}$ such that for some $p \in \mathcal{M}$, we have that
$E_p(f_1)\le E_p(f_2) \le\dots\le E_p(f_k)$.
\Ensure an index set of $\opt_{\succ}(\mathcal{K})$
\State $I \gets \{k\}$ \Comment{an index set of $\opt_{\succ}(\mathcal{K})$}
\For{$i=1:k-1$ }
\If{IsNotDominated2($i$)} 
\State{\indent $I\gets I \cup \{i\}$}\Comment{$f_i$ is maximal}
\EndIf
\EndFor
\State \Return $I$
\State{\textbf{where} IsNotDominated2($i$)}
\State{\indent \textbf{for} $j\in \{k,k-1,\dots,i+1\}$} 
\State{\indent \indent \textbf{if} $\underline{E}(f_j-f_i) > 0$  \textbf{then return False}} \Comment{$f_i$ is dominated by $f_j$}
\State{\indent\Return \True}

\end{algorithmic}
\end{algorithm}

Even though we have to do extra work to sort gambles at the beginning, we do not have to make as many comparisons in \cref{alg:FindMax2} as in \cref{alg:FindMax:matt}.
In particular, in the case that the set $\mathcal{K}$ has one maximal gamble, \cref{alg:FindMax2} only needs to solve $k-1$ linear programs. On the other hand, if all gambles are maximal, \cref{alg:FindMax2} needs to evaluate $\frac{k(k-1)}{2}$ linear programs. In both cases, this is only half of the number of comparisons of \cref{alg:FindMax:matt}.

Instead of solving multiple linear programs, \citet{2017:Jansen:Augustin:Schollmeyer} suggest to solve just a single linear program per gamble in $\mathcal{K}$:
\begin{align}
\label{P0:a}\tag{P0a}
\linprogref{P0} \quad 
\max \quad   & \sum_{j =1}^{k}\sum_{\omega \in \Omega}p_j(\omega)  \\
\label{P0:b}\tag{P0b}
\text{subject to} \quad & \forall j=1,\dots,k \colon  \sum_{\omega \in \Omega}p_j(\omega) \leq 1\\
\label{P0:c}\tag{P0c}
&  \forall j=1,\dots,k,\ \forall g_i\in\dom\underline{P} \colon \sum_{\omega \in \Omega} (g_i(\omega)-\underline{P}(g_i))p_j(\omega) \geq 0\\
\label{P0:d}\tag{P0d}
& \forall j=1,\dots,k \colon \sum_{\omega \in \Omega}\left(f(\omega)-f_j(\omega)\right)p_j(\omega) \geq 0\\
\label{P0:e}\tag{P0e}
\text{where} \quad & \forall j=1,\dots,k,\ \forall \omega \colon p_j(\omega)\geq 0.
\end{align}
If the optimal value of \linprogref{P0} is equal to $k$, then $f$ is a maximal gamble in $\mathcal{K}$ \citep{2017:Jansen:Augustin:Schollmeyer}. 
Therefore, to determine those $k$ gambles, we solve only $k$ linear programs (see \cref{alg:OneLP1}). However, the size of linear program is much bigger as it has $k(3+n)$ constraints and $mk$ variables.

Note that if we modify the constraint \cref{P0:b} to the following equality:
\begin{equation}\tag{P0b'}\label{P0:b'}
\forall j=1,\dots,k \colon  \sum_{\omega \in \Omega}p_j(\omega) = 1,
\end{equation}
then every feasible solution of \linprogref{P0'}:
\begin{equation}
\linprogref{P0'} \quad \max 0 \text{ subject to \cref{P0:b',P0:c,P0:d,P0:e}},
\end{equation} 
is also an optimal solution of \linprogref{P0}.
Therefore, if \linprogref{P0'} has a feasible solution, then $f$ is a maximal gamble in $\mathcal{K}$. In our simulation study, we solve \linprogref{P0'} because it is in a more suitable format for the primal-dual method, as it needs fewer artificial variables.

\begin{algorithm}
	\caption{Find the set of maximal gambles in $\mathcal{K}$}\label{alg:OneLP1}
	\begin{algorithmic}[1]
		\Require a set of $k$ gambles $\mathcal{K} = \{ f_1,\dots, f_k\}$
		\Ensure an index set of $\opt_{\succ}(\mathcal{K})$
		\State $I \gets \emptyset$ \Comment{an index set of $\opt_{\succ}(\mathcal{K})$}
			\For{$i = 1\colon k$}
				\If{\linprogref{P0'} with respect to $\mathcal{K}$ and $f_i$ has a feasible solution}
					\State{$I\gets I \cup \{i\}$}\Comment{$f_i$ is maximal}
					\EndIf
					\EndFor
		\State \Return $I$
	\end{algorithmic}
\end{algorithm}

Remember that if a gamble is not maximal in a given iteration, then we can exclude it from all further iterations. We use this idea to improve algorithm 3 and present it in \cref{alg:OneLP2}.
\begin{algorithm}
	\caption{Find the set of maximal gambles in $\mathcal{K}$}\label{alg:OneLP2}
	\begin{algorithmic}[1]
		\Require a set of $k$ gambles $\mathcal{K} = \{ f_1,\dots, f_k\}$
		\Ensure an index set of $\opt_{\succ}(\mathcal{K})$
		\State $I \gets \emptyset$ \Comment{an index set of $\opt_{\succ}(\mathcal{K})$}
		\For{$i = 1\colon k$}
		\If{IsNotDominated4$(I,i)$}
		\State{\indent $I\gets I \cup \{i\}$}\Comment{$f_i$ is maximal}
		\EndIf
		\EndFor
		\State \Return $I$
		\State{\textbf{where} IsNotDominated4($I,i$)}
		\State{\indent $\mathcal{G} = \{ f_j\in \mathcal{K}\colon j\in I \cup \{i+1,\dots,k\}\}$} 
		\State{\indent \textbf{if} \linprogref{P0'} with respect to $\mathcal{G}$ and $f_i$ has a feasible solution \textbf{then}}
		 \State{\indent\indent \textbf{return True}} 
	\end{algorithmic}
\end{algorithm}

\Cref{alg:FindMax:matt,alg:FindMax2,alg:OneLP1,alg:OneLP2} will be benchmarked later in \cref{sec:benchmark}.

\subsection{Algorithm for finding interval dominant gambles}

As we mentioned before, every maximal gamble in $\mathcal{K}$ is also interval dominant. Therefore, before running each algorithm, we can eliminate some non-maximal gambles in $\mathcal{K}$ by finding $\opt_{\sqsupset}(\mathcal{K})$. To check whether a gamble $f$ is interval dominant in $\mathcal{K}$, we first calculate $\max_{g\in \mathcal{K}}\underline{E}(g)$. Then $f$ is interval dominant if 
\begin{equation}\label{eq:check_id}
\overline{E}(f) \geq \max_{g\in \mathcal{K}}\underline{E}(g).
\end{equation}
Overall, to handle $k$ gambles, we have to solve $2k-1$ linear programs \citep[p.~337]{2014:troffaes:itip:computation}. This algorithm for finding interval dominant gambles in $\mathcal{K}$ is summarized in \cref{alg:FindID}.
\begin{algorithm}
\caption{Find the set of interval dominant gambles in $\mathcal{K}$}\label{alg:FindID}
\begin{algorithmic}[1]
\Require a set of $k$ gambles $\mathcal{K} = \{f_1,\dots, f_k \}$
\Ensure an index set of  $\opt_{\sqsupset}(\mathcal{K})$
\For{$j \in \{1,2,\dots,k\}$}
\State $e_j\gets\underline{E}(f_j)$
\EndFor
\State $\ell \gets \argmax_{j=1}^{k} e_j$
\State $I \gets \{\ell \}$ \Comment{an index set of  $\opt_{\sqsupset}(\mathcal{K})$}
\For{$i \in \{1,2,\dots,k\}\setminus \{\ell \}$}
\If {$\overline{E}(f_i) \geq e_\ell$}
\State{ $I \gets I \cup \{i\}$}\Comment{$f_i$ is interval dominant}
\EndIf
\EndFor
\State \Return $I$
\end{algorithmic}
\end{algorithm}
So, in \cref{sec:benchmark}, in addition to benchmarking those three algorithms for finding $\opt_{\succ}(\mathcal{K})$, we will also run \cref{alg:FindID} to helping those three algorithms to identify maximal gambles. Specifically, we will run \cref{alg:FindID} at the beginning to eliminate non-maximal gambles in $\mathcal{K}$, and then run those three algorithms on $\opt_{\sqsupset}(\mathcal{K})$.

\subsection{Fast evaluation of natural extensions inside algorithms}
\label{sec:fastnatext}

We can also speed up the process of evaluating the natural extension through \linprogref{P1} or \linprogref{D1}. To do so, we exploit the fact that we only need to find the sign of $\underline{E}(g-f)$, and not its exact value, to verify whether $f$ is dominated by $g$ or not.

As we minimize the objective function in \linprogref{P1}, the optimal value of \linprogref{P1} is less or equal to other feasible objective values. So, we can stop as soon as we find a feasible solution that achieves a negative objective value, because then we know that the optimal value of \linprogref{P1} is negative. In this case, $f$ is not maximal in $\mathcal{K}$.

Similarly, as we maximize the objective function in \linprogref{D1}, the optimal value of \linprogref{D1} is larger or equal to other feasible objective values. Consequently, we can stop as soon as we find a feasible solution that achieves a positive objective value.
In this case, $f$ is not dominated by $g$, so we must continue to compare $f$ to other gambles in $\mathcal{K}$.

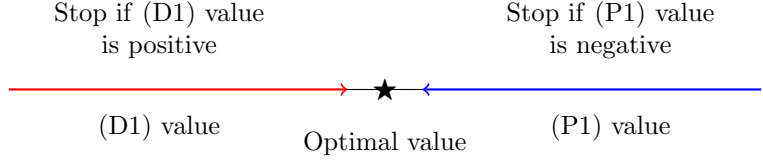
\begin{figure}
\centering
\begin{tikzpicture}
\draw (2,0.8) node[align=center] {Stop if \linprogref{D1} value \\  is positive};
\draw (8,0.8) node[align=center] {Stop if \linprogref{P1} value \\  is negative}; 
\draw (5,-1) node[anchor=south]{Optimal value};  
\draw (5,-0.01) node{$\bigstar$};
\draw (0,0)--(10,0);
\draw [thick,red,->] (0,0)--(4.5,0);
\draw (2,-0.8) node[anchor=south]{ \linprogref{D1}  value}; 
\draw [thick,blue,<-] (5.5,0)--(10,0);
\draw (8,-0.8) node[anchor=south]{ \linprogref{P1}  value}; 
\end{tikzpicture}
\caption{Early stopping criterion}\label{fig:1}
\end{figure}

These extra stopping criteria are illustrated in \cref{fig:1}.

We can solve \linprogref{P1} and \linprogref{D1} by many linear programming methods, for example, the simplex, the affine scaling or the primal-dual interior point methods. However, we only focus on the primal-dual method as this method solves both primal and dual problems simultaneously, and can therefore exploit both stopping criteria simultaneously. We also know from earlier work that the primal-dual interior point method is particularly suitable for working with lower previsions \citep{2018:Nakharutai:Troffaes:Caiado}. Finally, the primal-dual interior point method is widely regarded as one of the best general purpose linear programming methods \cite[\S 10.2]{2009:Griva:Nash:Sofer}.

Note that, in practice, the primal-dual method can start with an arbitrary point and then generates a sequence of (not necessarily feasible) points that converges to an optimal feasible solution \citep[\S 7.3]{1993:Fang:Puthenpura}. On the other hand, given initial feasible points, the method will generate a sequence of feasible points converging to an optimal solution \citep[\S 7.3]{1993:Fang:Puthenpura}. Therefore, 
to apply the extra stopping criterion, we first have to find initial feasible points for \linprogref{P1} and \linprogref{D1}.
Fortunately, there is an efficient way to obtain initial feasible points for the linear programs \linprogref{P1} and \linprogref{D1}.

For \linprogref{P1}, we can apply the first phase of the two-phase method to obtain an initial feasible probability mass function $p(\omega)$ \citep[\S 4.2]{2018:Nakharutai:Troffaes:Caiado}. This technique is usually used for obtaining interior feasible points for the affine scaling method \citep[\S 7.1.2]{1993:Fang:Puthenpura}. 
Since the constraints of \linprogref{P1} do not change, once we find a feasible starting point, we can reuse it for other problems \linprogref{P1} with different objective functions. So we only need to do this once for any given lower prevision, and it is independent of the decision problem.

For \linprogref{D1}, we can very quickly calculate a feasible starting point without solving a linear program, using a result from \citet[Theorem~7]{2018:Nakharutai:Troffaes:Caiado}.

Unfortunately, there is no direct way to obtain feasible starting points for the linear programming problem \linprogref{P0'} (otherwise we would have immediately solved the problem).

\section{Benchmarking}\label{sec:benchmark}

\subsection{Generating sets of gambles for benchmarking}

As we mentioned before, we would like to generate a set of gambles $\mathcal{K}$ for benchmarking \cref{alg:FindMax:matt,alg:FindMax2,alg:OneLP1,alg:OneLP2} for finding $\opt_{\succ}(\mathcal{K})$ and \cref{alg:FindID} for finding $\opt_{\sqsupset}(\mathcal{K})$. Can we generate a set $\mathcal{K}$ such that $|\mathcal{K}|=k$, $|\opt_{\succ}(\mathcal{K})|=m$ and $|\opt_{\sqsupset}(\mathcal{K})| = n$ where $m\leq n \leq k$?

A naive idea is to first generate $\mathcal{K} = \{g\}$, so obviously, $\opt_{\succ}(\mathcal{K}) = \opt_{\sqsupset}(\mathcal{K})=\{g\}$. Next, we generate a gamble $h$ such that
\begin{align}\label{eq:h-maximal}
\opt_{\succ}(\mathcal{K}\cup \{ h\}) & = \opt_{\succ}(\mathcal{K}) \cup  \{ h\}\quad \&\quad \opt_{\sqsupset}(\mathcal{K}\cup \{ h\}) =\opt_{\sqsupset}(\mathcal{K})\cup \{ h\},\\ 
\intertext{and we add $h$ to $\mathcal{K}$. We repeat this process until we have $|\mathcal{K}|=|\opt_{\succ}(\mathcal{K})|=|\opt_{\sqsupset}(\mathcal{K})| = m$. After that, we generate a gamble $h$ that satisfies}\label{eq:h-idNotmax}
\opt_{\succ}(\mathcal{K}\cup \{ h\})& = \opt_{\succ}(\mathcal{K}) \quad \&\quad \opt_{\sqsupset}(\mathcal{K}\cup \{ h\}) =\opt_{\sqsupset}(\mathcal{K})\cup \{ h\}.\\ 
\intertext{Again, we add $h$ to $\mathcal{K}$ and repeat this process until we have $|\mathcal{K}|=|\opt_{\sqsupset}(\mathcal{K})| = n$. However, $|\opt_{\succ}(\mathcal{K})|= m$. Finally, we generate a gamble $h$ such that}\label{eq:h-Notid}
\opt_{\succ}(\mathcal{K}\cup \{ h\})& = \opt_{\succ}(\mathcal{K}) \quad \&\quad \opt_{\sqsupset}(\mathcal{K}\cup \{ h\}) =\opt_{\sqsupset}(\mathcal{K}),
\end{align}
we add $h$ to $\mathcal{K}$, and repeat this process until we have $|\mathcal{K}|=k$, $|\opt_{\succ}(\mathcal{K})|=m$ and $|\opt_{\sqsupset}(\mathcal{K})| = n$ as we want.

In practice, a randomly generated gamble $h$ may not easily satisfy \cref{eq:h-maximal}, \cref{eq:h-idNotmax} or \cref{eq:h-Notid}. We may need to sample many gambles until we satisfy the desired conditions, and therefore it may take a while to obtain the set $\mathcal{K}$ that we want.

Surprisingly, for any generated gamble $h$, we can modify $h$ by shifting it by $\alpha$, for some $\alpha \in \mathbb{R}$, so that a new gamble $h-\alpha$ meets any of the above requirements.
Next, we explain for what range of $\alpha$, the gamble $h-\alpha$ satisfies either \cref{eq:h-maximal}, \cref{eq:h-idNotmax} or \cref{eq:h-Notid}. Specifically, we identify for which values of $\alpha$ we have one of the following:
\begin{enumerate}
\item[(i)] $\opt_{\succ}(\mathcal{K}\cup \{ h-\alpha\}) = \opt_{\succ}(\mathcal{K}) \cup  \{ h-\alpha\}$, 
\item[(ii)] $\opt_{\succ}(\mathcal{K}\cup \{ h-\alpha\}) = \opt_{\succ}(\mathcal{K})$ and $\opt_{\sqsupset}(\mathcal{K}\cup \{ h-\alpha\}) =\opt_{\sqsupset}(\mathcal{K})\cup \{ h-\alpha\}$,
\item[(iii)] $\opt_{\sqsupset}(\mathcal{K}\cup \{ h-\alpha\}) =\opt_{\sqsupset}(\mathcal{K})$
\end{enumerate}

Let $\mathcal{K}$ be a set of gambles.
Given any gamble $h$, \cref{lem:upBoundAlpha} shows for which $\alpha$, $h-\alpha$ is a maximal gamble in $\mathcal{K} \cup \{h - \alpha\}$. 

\begin{lemma}\label{lem:upBoundAlpha}
Let $\mathcal{K}$ be a set of gambles and let $h$ be another gamble and $\alpha \in \mathbb{R}$. Then $h-\alpha$ is maximal in $\mathcal{K} \cup \{h - \alpha\}$ if and only if 
\begin{equation}
\min_{f \in \opt_{\succ}(\mathcal{K})} \overline{E}(h-f) \geq \alpha.
\end{equation}
\end{lemma}

\Cref{lem:upBoundAlpha} provides an upper bound on $\alpha$ for which $h-\alpha$ is maximal in $\mathcal{K} \cup \{h - \alpha\}$. However, if we set $\alpha$ too low, then $h-\alpha$ may dominate other maximal gambles in $\opt_{\succ}(\mathcal{K})$, that is, we risk having gambles $f$ for which $f \in \opt_{\succ}(\mathcal{K})$ but $f \notin \opt_{\succ}(\mathcal{K} \cup \{h - \alpha\})$. The following lemma tells us how to prevent this situation.

\begin{lemma}\label{lem:LowBoundAlpha}
Let $\mathcal{K}$ be a set of gambles and let $h$ be another gamble and $\alpha \in \mathbb{R}$. Then all maximal gambles in $\mathcal{K}$ are still maximal in $\mathcal{K}\cup \{h - \alpha\}$ if and only if
\begin{equation}
\max_{f \in \opt_{\succ}(\mathcal{K})}\underline{E}(h-f) \leq \alpha.
\end{equation}
\end{lemma}

\Cref{lem:LowBoundAlpha} provides a lower bound on $ \alpha$ such that $h -\alpha$ does not dominate any other maximal gambles in $\mathcal{K}\cup \{h - \alpha\}$.

Finally, by \cref{eq:check_id}, we know that $h-\alpha$ is interval dominant in $\mathcal{K}\cup \{h - \alpha\}$ if and only if
\begin{equation}
\overline{E}(h-\alpha) \geq \max_{f\in \mathcal{K}}\underline{E}(f).
\end{equation}
This is equivalent to 
\begin{equation}
\alpha \leq \overline{E}(h) - \max_{f \in \mathcal{K}}\underline{E}(f).
\end{equation}

The next lemma ensures that these bounds on $\alpha$ are always ordered in the same way:
\begin{lemma}\label{lem:max-min_h-f}
Let $\mathcal{K}$ be a set of gambles and let $h$ be another gamble. Then, the following holds:
\begin{equation}\label{eq:lem:max-min}
\max_{f \in \opt_{\succ}(\mathcal{K})}\underline{E}(h-f) \leq \min_{ f\in  \opt_{\succ}(\mathcal{K})} \overline{E}(h - f) \leq \overline{E}(h) - \max_{f \in \mathcal{K}}\underline{E}(f).
\end{equation}
\end{lemma}

\Cref{thm:range-alpha} brings everything together, and summarises for which ranges of $\alpha$ we have that $h-\alpha$ satisfies either \cref{eq:h-maximal}, \cref{eq:h-idNotmax} or \cref{eq:h-Notid}.

\begin{figure}
	\centering
	\begin{tikzpicture}[scale=0.8,transform shape]
	\draw (-1,0) -- (1,0)node[anchor=north]{$\displaystyle\max_{f \in \opt_{\succ}(\mathcal{K})}\underline{E}(h-f)$} --(6,0)node[anchor=north]{$\displaystyle\min_{ f\in  \opt_{\succ}(\mathcal{K})} \overline{E}(h - f)$} --(11,0)node[anchor=north]{$\overline{E}(h) -\displaystyle\max_{f \in \mathcal{K}}\underline{E}(f)$};
	\draw[->] (11,0) -- (13,0) node[anchor=west]{$\alpha$};
	\draw   (1.2,0.3)--(0.98,0.3);
	\draw (1,-0.1)--(1,0.3);
	\draw   (5.8,0.3)--(6.2,0.3);
	\draw (6,-0.1)--(6,0.3);
	\draw   (11,-0.1)--(11,0.3);
	\draw   (11.2,0.3) -- (10.8,0.3);
	\draw[->] (3.5,1.25) -- (3.5,0.25);
	\draw(3.5,2.5)node[anchor=north]{\small $\opt_{\succ}(\mathcal{K}\cup \{ h-\alpha\}) = \opt_{\succ}(\mathcal{K}) \cup  \{ h-\alpha\},$};
	\draw(3.5,2)node[below =20pt of {(3.5,2.5)},anchor=north]{\small $\opt_{\sqsupset}(\mathcal{K}\cup \{ h-\alpha\}) =\opt_{\sqsupset}(\mathcal{K})\cup \{ h-\alpha\}$};
	\draw[->] (8.5,2.8) -- (8.5,0.25);
	\draw(7.7,4.2)node[anchor=north]{\small $\opt_{\succ}(\mathcal{K}\cup \{ h-\alpha\}) = \opt_{\succ}(\mathcal{K}),$};  
	\draw(8.5,3.5)node[anchor=north]{\small $\opt_{\sqsupset}(\mathcal{K}\cup \{ h-\alpha\}) =\opt_{\sqsupset}(\mathcal{K})\cup \{ h-\alpha\}$};
	\draw[->] (12.4,1.2) -- (12.4,0.25);
	\draw(11.5,2.5)node[anchor=north]{\small $\opt_{\succ}(\mathcal{K}\cup \{ h-\alpha\}) = \opt_{\succ}(\mathcal{K}),$};  
	\draw(11.5,1.8)node[anchor=north]{\small $\opt_{\sqsupset}(\mathcal{K}\cup \{ h-\alpha\}) =\opt_{\sqsupset}(\mathcal{K})$};    
	\end{tikzpicture}
	\caption{Ranges for $\alpha$ such that $h-\alpha$ satisfies either of the situations described in \cref{thm:range-alpha}.}\label{fig:alpha}
\end{figure}
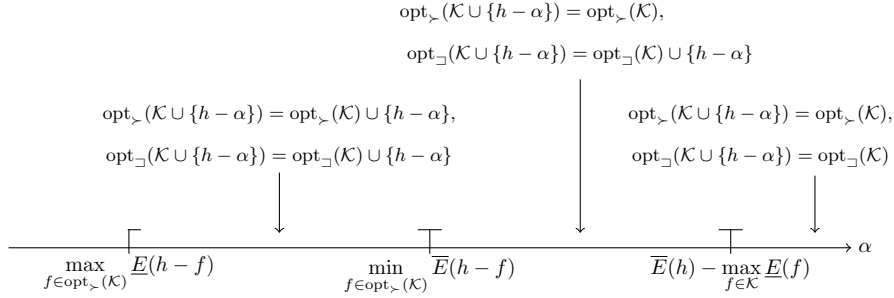

\begin{theorem}\label{thm:range-alpha}
Let $\mathcal{K}$ be a set of gambles and let $h$ be another gamble and $\alpha \in \mathbb{R}$.
\begin{enumerate}
	\item If we choose
	\begin{equation}\label{eq:rangeAlpha}
\max_{f \in \opt_{\succ}(\mathcal{K})}\underline{E}(h-f)\leq \alpha\leq \min_{f \in \opt_{\succ}(\mathcal{K})} \overline{E}(h-f)
	\end{equation}
	then
\begin{align}
\opt_{\succ}(\mathcal{K}\cup \{ h-\alpha\})
&=\opt_{\succ}(\mathcal{K}) \cup  \{ h-\alpha\} \\
\opt_{\sqsupset}(\mathcal{K}\cup \{ h-\alpha\}) 
&=\opt_{\sqsupset}(\mathcal{K})\cup \{ h-\alpha\}.
\end{align}
	\item If we choose
	\begin{equation}\label{eq:alpha-condi2}
\min_{ f\in  \opt_{\succ}(\mathcal{K})} \overline{E}(h - f) < 	\alpha \leq \overline{E}(h) - \max_{f \in \mathcal{K}}\underline{E}(f),
	\end{equation}
	then
\begin{align}
\opt_{\succ}(\mathcal{K}\cup \{ h-\alpha\})
&=\opt_{\succ}(\mathcal{K}) \\
\opt_{\sqsupset}(\mathcal{K}\cup \{ h-\alpha\}) 
&=\opt_{\sqsupset}(\mathcal{K})\cup \{ h-\alpha\}.
\end{align}
	\item If we choose 
	\begin{equation}
	\alpha > \overline{E}(h) -\max_{f \in \mathcal{K}}\underline{E}(f),
	\end{equation}
	then
\begin{align}
\opt_{\succ}(\mathcal{K}\cup \{ h-\alpha\})
&= \opt_{\succ}(\mathcal{K}) \\
\opt_{\sqsupset}(\mathcal{K}\cup \{ h-\alpha\}) &=\opt_{\sqsupset}(\mathcal{K}).
\end{align}
\end{enumerate}	
\end{theorem}
\Cref{fig:alpha} illustrates \cref{thm:range-alpha}.


From \cref{thm:range-alpha}, we obtain an algorithm for generating a set $\mathcal{K}$ of $k$ gambles such that $|\opt_{\succ}(\mathcal{K})| = m$ and $|\opt_{\sqsupset}(\mathcal{K})| = n$ where $m\leq n \leq k$.
First, we generate a set $\mathcal{K}$ consisting of $m$ gambles that are all maximal,
so $|\opt_{\succ}(\mathcal{K})|=|\opt_{\sqsupset}(\mathcal{K})|= |\mathcal{K}|=m$.
Next, we add $n-m$ further gambles to $\mathcal{K}$, where these gambles are interval dominant but not maximal,
so $|\opt_{\sqsupset}(\mathcal{K})|= |\mathcal{K}|=n$ but it remains that $|\opt_{\succ}(\mathcal{K})|=m$.
Finally, we add $k-n$ further gambles to $\mathcal{K}$, where these gambles are not interval dominant.
We then end up with a set of gambles which has a precisely given number of maximal and interval dominant gambles, as required.
For the full algorithm, see \cref{alg:generate-set}.

\begin{algorithm}
\caption{Generate a set of $k$ gambles $\mathcal{K}$ such that $|\opt_{\succ}(\mathcal{K})| = m$ and $|\opt_{\sqsupset}(\mathcal{K})| = n$ where $m\leq n \leq k$}\label{alg:generate-set}
\begin{algorithmic}[1]
\Require (a) Numbers $m$, $n$, $k$ where $m\leq n \leq k$, and (b)
a sequence of $k$ gambles $h_1$, \dots, $h_k$ such that $\overline{E}(h_i-h_j)<\overline{E}(h_i)-\underline{E}(h_j)$ for all $i$, $j\in\{1,\dots,k\}$
\Ensure a set of $k$ gambles $\mathcal{K}$ such that such that $|\opt_{\succ}(\mathcal{K})| = m$ and $|\opt_{\sqsupset}(\mathcal{K})| = n$ where $m\leq n \leq k$
\State{$\mathcal{K} \gets \{h_1\}$} 
\For{$i = 2:m$}  \Comment{$m$ maximal and interval dominant}
\State Choose $\alpha$ such that $ \displaystyle{\max_{f \in \opt_{\succ}(\mathcal{K})}}\underline{E}(h_i-f)\leq \alpha\leq \min_{f \in \opt_{\succ}(\mathcal{K})} \overline{E}(h_i-f)$
\State $\mathcal{K}\gets  \mathcal{K} \cup \{h_i-\alpha\}$
\EndFor
\For{$i = m+1:n$} \Comment{\small $n-m$ interval dominant but not maximal}\normalsize
\State Choose $\alpha$ such that $\displaystyle{\min_{f \in \opt_{\succ}(\mathcal{K})} \overline{E}(h_i-f)< \alpha\leq \overline{E}(h_i) - \max_{f\in \mathcal{K}}\underline{E}(f)}$
\State $\mathcal{K}\gets  \mathcal{K} \cup \{h_i - \alpha\}$ 
\EndFor
\For{$i = n+1:k$} \Comment{\small $k-n$ not interval dominant}\normalsize
\State Choose $\alpha$ such that $\alpha  >\overline{E}(h_i) - \displaystyle\max_{f\in \mathcal{K}}\underline{E}(f)$
\State $\mathcal{K}\gets  \mathcal{K} \cup \{h_i - \alpha\}$ 
\EndFor
\State \Return $\mathcal{K}$
\end{algorithmic}
\end{algorithm}

Note that if $\alpha$ in the third loop is much larger than $\overline{E}(h)-\max_{f\in \mathcal{K}}\underline{E}(f)$, then $h - \alpha$ can be more easily dominated.
Therefore, if we  want these non-maximal gambles to be difficult to detect, then we should set $\alpha$ to be only slightly larger than $\overline{E}(h)-\max_{f\in \mathcal{K}}\underline{E}(f)$.
In our simulation study, we will choose $\alpha$ in the third loop to be $\overline{E}(h)-\max_{f\in \mathcal{K}}\underline{E}(f) + \epsilon$ where $\epsilon$ is uniformly sampled from $(0,1)$.

Also note that we require the following condition for all $i$ and $j$:
\begin{equation}\label{eq:strictineqfordifference}
\overline{E}(h_i-h_j)<\overline{E}(h_i)-\underline{E}(h_j).
\end{equation}
This ensures that there exists an $\alpha$ satisfying the strict inequality in \cref{eq:alpha-condi2}, because in that case, the left hand side of \cref{eq:alpha-condi2} will be strictly less than the right hand side of \cref{eq:alpha-condi2} (see \cref{eq:thelastone} in the appendix; the inequality there will be a strict inequality under the assumed condition). In this way, there is always an $\alpha$ for which $h-\alpha$ is not maximal but still interval dominant.
\Cref{eq:strictineqfordifference} requires that $\underline{E}$
is non-linear (i.e. genuinely imprecise),
and that the gambles $h_i$ are non-constant and linearly independent. For example, if for each $\omega$, we sample $h_i(\omega)$ uniformly from $[0,1]$,
then this requirement is practically always satisfied, and if not, we can simply resample $h_i$
until it is.

Also note that \cref{alg:generate-set} has to evaluate many natural extensions for each new generated gamble. This is required to ensure that the numbers of maximal gambles and interval dominant gambles in the sets are exactly as specified.
Although this gives us very precise control
over the range of decision problems that are generated, an obvious
downside is that evaluating all these natural extensions requires a
huge computational effort. For this reason, we had to limit ourselves
to $|\Omega| \leq 2^6$ and $|\mathcal{K}|\leq 2^8$.

\subsection{Benchmarking results}

\begin{figure}[p]
\centering
\begin{tikzpicture}
\filldraw[draw=black, fill=gray!20] (2,0) -- (8,4.5) --(8,0);
\draw (2,-1) --(2,5.5)node[anchor=east]{$m$};
\draw (2,4.5) --(9.5,4.5);
\draw [->](2,5.5) -- (2,6);
\draw (6,3) --(8,3);
\draw (4,1.5) -- (8,1.5);
\draw (2,0)--(8,4.5);
\draw (8,0)-- (8,1.5)node[circle,draw=black,fill=white!80!black,minimum size=20]{$g$} --(8,3)node[circle,draw=black,fill=white!80!black,minimum size=20]{$i$}--(8,4.5)node[circle,draw=black,fill=white!80!black,minimum size=20]{$j$} --(8,5.5);
\draw (8.5,5.2)node[anchor=west]{$k$};
\draw(8,5) -- (8,6);

\draw (6,0)-- (6,1.5)node[circle,draw=black,fill=white!80!black,minimum size=20]{$f$} --(6,3)node[circle,draw=black,fill=white!80!black,minimum size=20]{$h$};

\draw (4,0)-- (4,1.5)node[circle,draw=black,fill=white!80!black,minimum size=20]{$e$};

\draw (1,0) -- (2,0)node[circle,draw=black,fill=white!80!black,minimum size=20]{$a$} --(4,0)node[circle,draw=black,fill=white!80!black,minimum size=20]{$b$}--(6,0)node[circle,draw=black,fill=white!80!black,minimum size=20]{$c$}--(8,0)node[circle,draw=black,fill=white!80!black,minimum size=20]{$d$} --(10,0)node[anchor=north]{$n$};
\draw [->](10,0) -- (11,0);
\end{tikzpicture}
\caption{The area of $m\leq n \leq k$ and 10 options label the different $m$ and $n$ that we consider in the simulation (see \cref{tab:10cases})}\label{fig:mnk}
\end{figure}
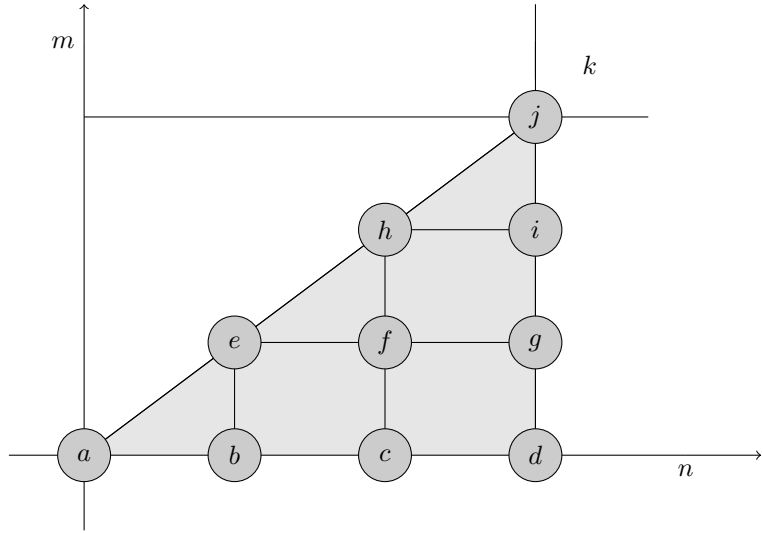

To benchmark those \cref{alg:FindMax:matt,alg:FindMax2,alg:OneLP1,alg:OneLP2,alg:FindID} from \cref{sec:algorithm}, in this section, we generate random sets of gambles. We consider the case that $|\Omega| = 2^2$ and $|\Omega| = 2^6$ and the number of gambles in $\mathcal{K}$ where $k = 2^j$ for $j \in \{4,6,8\}$. For each case, random sets of gambles $\mathcal{K}$ are generated as follows.  
\begin{enumerate}
\item We first generate a lower prevision $\underline{P}$ on a finite domain, that avoids sure loss. To do so, we use \citep[algorithm 2]{2018:Nakharutai:Troffaes:Caiado} with $2^4$ coherent previsions to generate a lower prevision $\underline{E}$ on the set of all gambles, that avoids sure loss. Next, we use \citep[stages 1 and 2 in algorithm 4]{2018:Nakharutai:Troffaes:Caiado} to restrict $\underline{E}$ to a lower prevision $\underline{P}$ that avoids sure loss, with a given finite size of domain. In this simulation we consider
$|\dom\underline{P}| = 2^i$ for $i \in \{2,4,6\}$.
This ensures that the generated lower prevision has no specific structural properties (such as 2-monotonicity \citep[Chapter~6]{2014:troffaes:decooman::lower:previsions}).
\item We generate $k$ gambles $h_1$, \dots, $h_k$ as follows. For each $\omega$ and $i$, we sample $h_i(\omega)$ uniformly from $[0,1]$ and check whether they satisfy \cref{eq:strictineqfordifference}.
\item We use \cref{alg:generate-set} to generate random sets $\mathcal{K}$ such that $|\mathcal{K}|=k$, $|\opt_{\succ}(\mathcal{K})| = m$ and $|\opt_{\sqsupset}(\mathcal{K})| = n$ where $m \leq n \leq k$, where we use the previously generated $\underline{P}$ to evaluate  $\underline{E}$ and $\overline{E}$.
Note that in \cref{alg:generate-set}, we choose $\alpha$ in the first loop as follows: we sample $\delta$ uniformly from $(0,1)$, and set
 \begin{equation}
 \alpha\coloneqq\delta\max_{f \in \opt_{\succ}(\mathcal{K})}\underline{E}(h_i-f)+(1-\delta)\min_{f \in \opt_{\succ}(\mathcal{K})} \overline{E}(h_i-f).  
 \end{equation}
 For $\alpha$ in the second loop, we choose it as follow: sample $\delta$ uniformly from $(0,1)$, and set
  \begin{equation}
 \alpha\coloneqq\delta\min_{f \in \opt_{\succ}(\mathcal{K})} \overline{E}(h_i-f)+(1-\delta)\left(\overline{E}(h_i) - \max_{f\in \mathcal{K}}\underline{E}(f)\right),
 \end{equation}
 and in the last loop, we set $\alpha\coloneqq\overline{E}(h_i) - \max_{f\in \mathcal{K}}\underline{E}(f) + \epsilon$, where we sample $\epsilon$ uniformly from $(0,1)$. 
\end{enumerate}  
In the simulation, we would like to cover a range of possible options of $m$, $n$, and $k$ that satisfy $m\leq n \leq k$. For each different size of $\mathcal{K}$, we consider $10$ options that vary the number of maximal gambles $m$ and the number of interval dominant gambles $n$ in $\mathcal{K}$  which are illustrated in \cref{fig:mnk} and \cref{tab:10cases}.

{
\setlength{\arrayrulewidth}{0.5mm}
\setlength{\tabcolsep}{12pt}
\renewcommand{\arraystretch}{1.3}
\begin{table}[p]
	\centering
	\begin{tabular}{|c|c|cV{1.5}c|cV{1.5}c|c|}
		\hline
		\multirow{2}{*}{Options} & \multicolumn{2}{lV{1.5}}{$|\mathcal{K}| = 2^4$} & \multicolumn{2}{lV{1.5}}{$|\mathcal{K}| = 2^6$} & \multicolumn{2}{l|}{$|\mathcal{K}| = 2^8$} \\ \cline{2-7} 
		&  $m$  &  $n$  &     $m$  &  $n$     &   $m$  &  $n$ \\ \hline
		\circled{a} &  1 & 1 &  1   & 1 &  1 & 1\\ \hline
		\circled{b} &  1 & 5  &  1 &  21 &   1 & 85\\ \hline
		\circled{c} &  1 & 11 & 1  &  42 &  1  &170 \\ \hline
		\circled{d} &  1 & 16 &  1 & 64  &  1  &256 \\ \hline
		\circled{e} &  5 & 5  &  21 &  21 &   85 &85 \\ \hline
		\circled{f} &  5 & 11 &  21 &  42 &  85  &170 \\ \hline
		\circled{g} &  5 & 16 &  21 & 64  &   85 & 256\\ \hline
		\circled{h} &  11 & 11 &  42 &  42 &  170  &170 \\ \hline
		\circled{i} &  11 & 16 &  42 &  64 &   170 & 256\\ \hline
		\circled{j}&  16 & 16  &  64 & 64  &  256  & 256\\ \hline
	\end{tabular}
	\caption{Table of points that indicate different sizes of set $\mathcal{K}$ with vary the number of maximal gambles $m$ and the number of interval dominant gambles $n$ in $\mathcal{K}$} \label{tab:10cases}
\end{table}
}

These 10 options can be grouped as follows. Options a to d represent the cases where $m=1$ while we increase $n$ from $1$ to $k$. Options d, g, i, and j represent the cases where $n=k$ while we increase $m$ from $1$ to $k$. Options a, e, h, and j represent the cases where $m=n$ while we increase them jointly from $1$ to $k$. Option f represents a case where $m<n<k$.

We then apply \cref{alg:FindMax:matt,alg:FindMax2,alg:OneLP1,alg:OneLP2,alg:FindID} on each generated set of gambles $\mathcal{K}$. For \cref{alg:FindMax:matt,alg:FindMax2,alg:FindID}, we solve linear programs for evaluating  upper and lower natural extensions by the improved primal-dual interior point method, including all of the improvements discussed in \cref{sec:fastnatext}, i.e. feasible starting points and early stopping criteria. As the
MATLAB implementation of the primal-dual method cannot be easily modified, we had to write our own implementation of the improved primal-dual interior point method in MATLAB (R2018a) \citep{MATLAB:2018}, based on the implementation used in \citet{2018:Nakharutai:Troffaes:Caiado}. For \cref{alg:OneLP1,alg:OneLP2}, we simply solve linear programs by the standard primal-dual method as the improvements cannot be applied. For fair comparison, we also used our own implementation of the standard primal-dual method here, rather than using the one that is available in MATLAB (R2018a).

\begin{figure}
	\centering
	\setlength{\tabcolsep}{2pt}
	\newcolumntype{C}{>{\centering\arraybackslash} m{0.48\linewidth} }

	\begin{tabular}{m{0.5em}CC}
		&
	$|\Omega| = 2^2$
		&
		$|\Omega| = 2^6$
		\\
		\rotatebox[origin=l]{90}{$|\mathcal{K}| = 2^4$}
		&
		\includegraphics[width=\hsize, trim={2cm 0 2cm 0},clip]{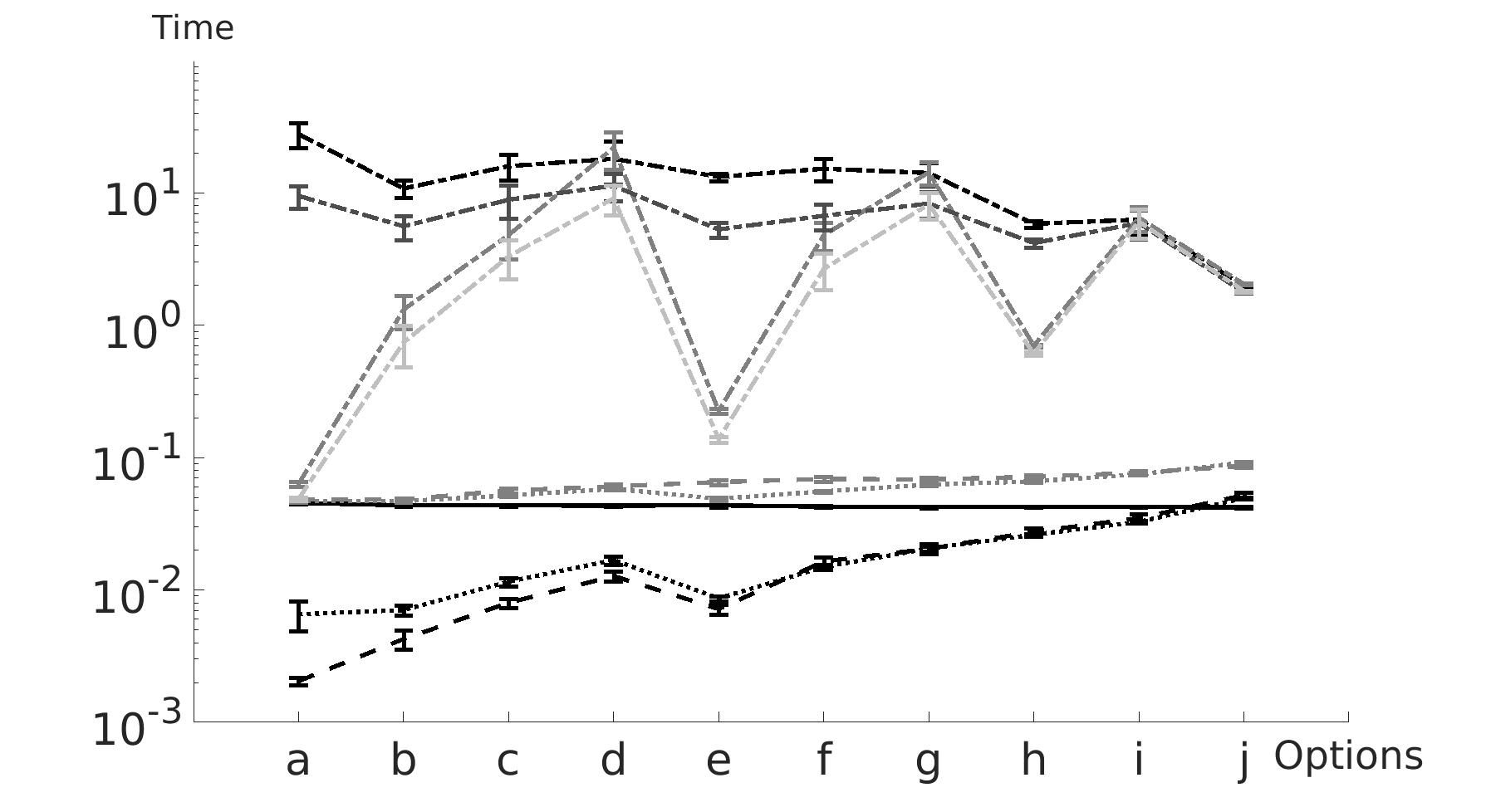}
		&
		\includegraphics[width=\hsize, trim={2cm 0 2cm 0},clip]{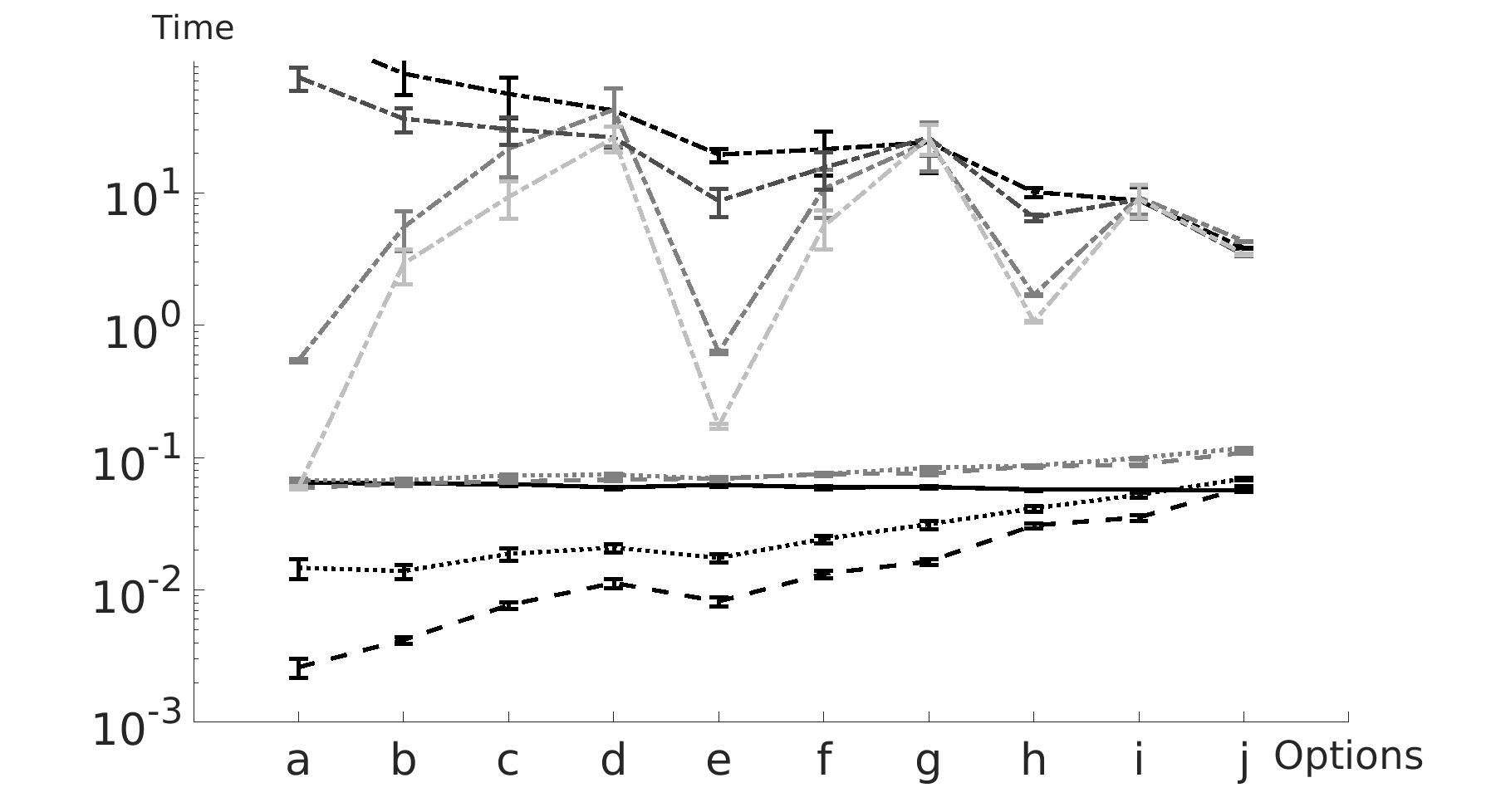}
		\\
		\rotatebox[origin=l]{90}{$|\mathcal{K}| = 2^6$}
		&
		\includegraphics[width=\hsize, trim={2cm 0 2cm 0},clip]{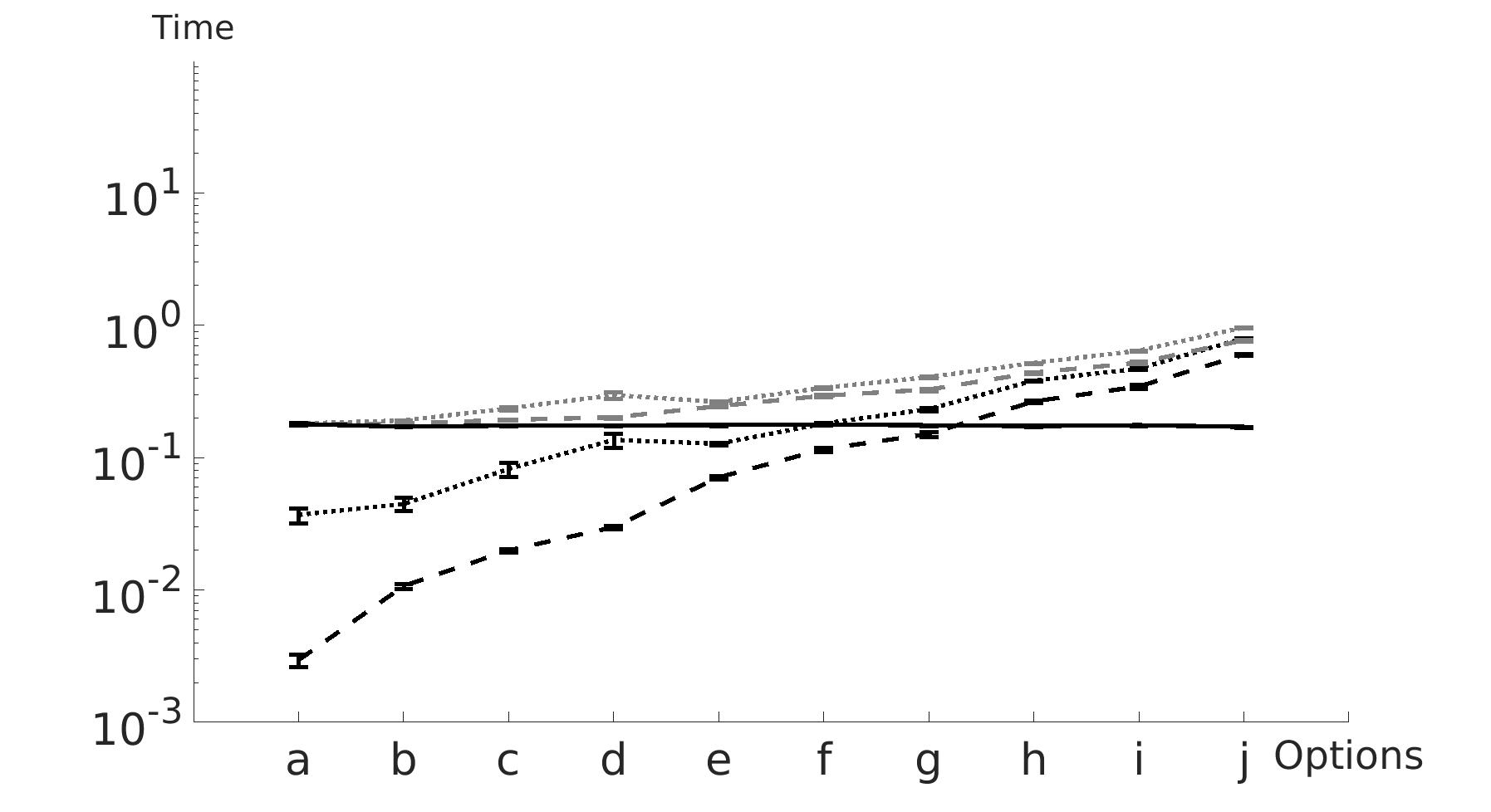}
		&
		\includegraphics[width=\hsize, trim={2cm 0 2cm 0},clip]{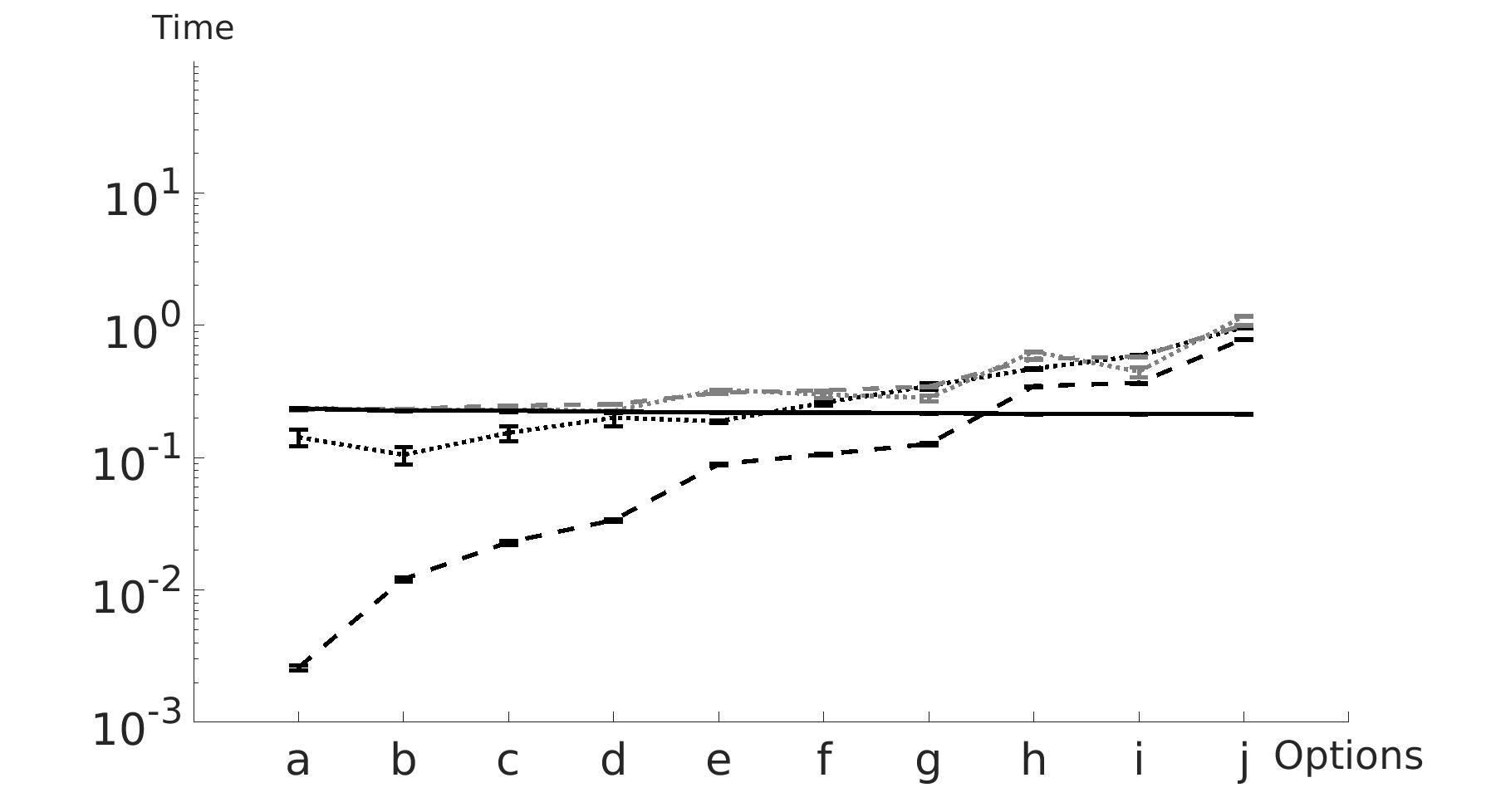}
		\\
		\rotatebox[origin=l]{90}{$|\mathcal{K}| = 2^8$}
		&
		\includegraphics[width=\hsize, trim={2cm 0 2cm 0},clip]{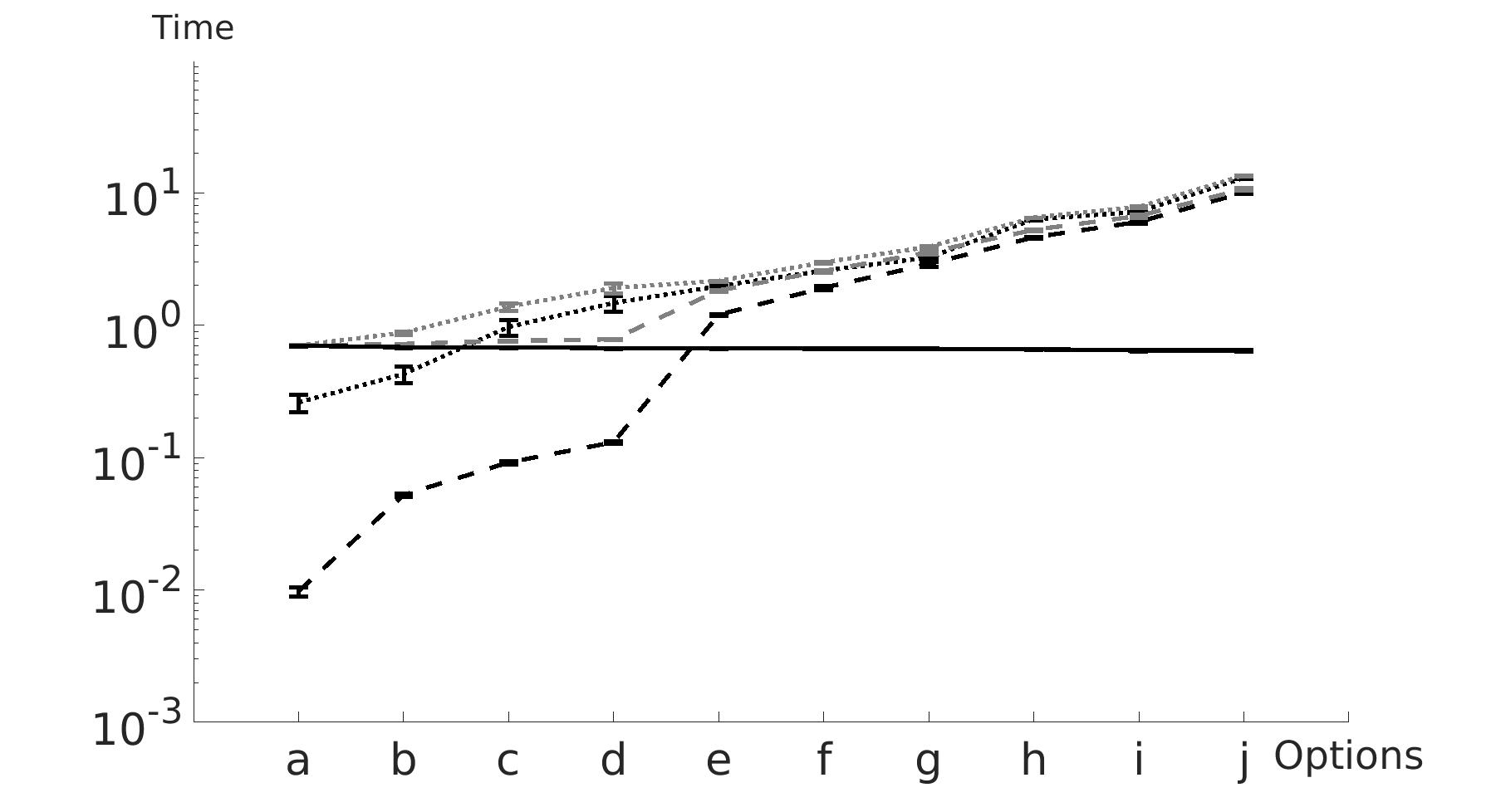}
		&
		\includegraphics[width=\hsize, trim={2cm 0 2cm 0},clip]{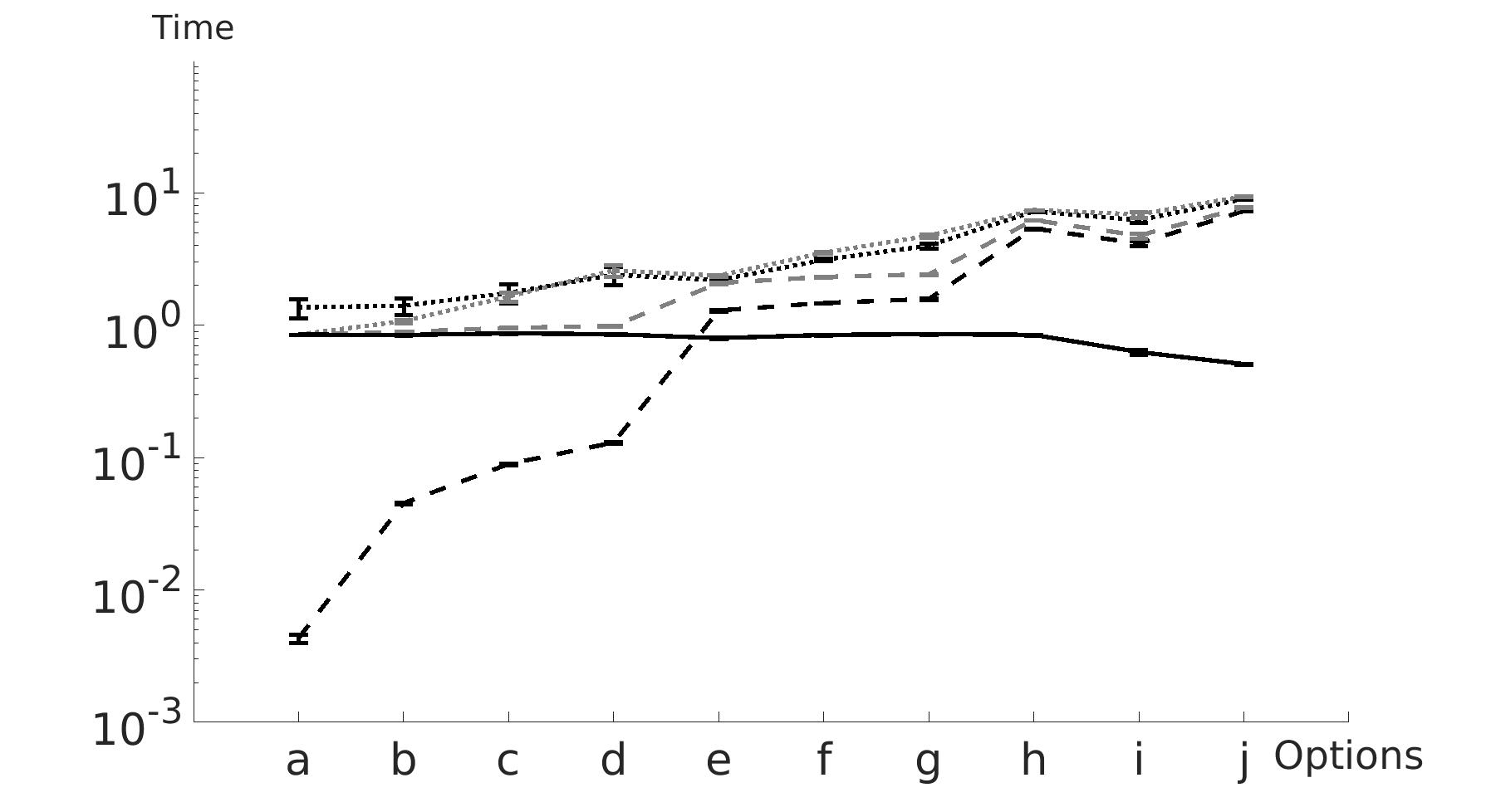}
	\end{tabular}
	\includegraphics[width=0.9\linewidth]{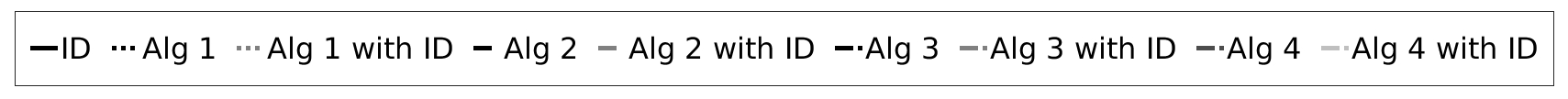}
	\caption{Comparison plots of the average computational time for \cref{alg:FindMax:matt,alg:FindMax2,alg:OneLP1,alg:OneLP2} for finding maximal gambles and for \cref{alg:FindID} for finding interval dominant gambles. The number of outcomes in left column is $2^2$ and $2^6$ in the right column.  Each row represents a different number of gambles with vary options of the numbers of maximal gambles and interval dominant gambles in the set (see \cref{tab:10cases} for each option). The labels indicate algorithms with and without  \cref{alg:FindID}. We fix $|\dom \underline{P}| = 2^4$.}
	\label{fig:plot1}
\end{figure}

To investigate whether interval dominance is helpful for finding maximal gambles, we also run each algorithm with and without \cref{alg:FindID}. Specifically, we run \cref{alg:FindMax:matt,alg:FindMax2,alg:OneLP1,alg:OneLP2} as such,
but additionally we also run \cref{alg:FindID} to obtain a set of interval dominant gambles, and then again run each of \cref{alg:FindMax:matt,alg:FindMax2,alg:OneLP1,alg:OneLP2} on just the resulting set of interval dominant gambles. In all cases, we measure the total computational time taken, i.e. including the time taken on \cref{alg:FindID} when applicable. Note that computational time to run \cref{alg:FindMax2} includes the time for sorting gambles from the lowest to the highest expectation as in \cref{eq:E_p(f_j-f_i)}. To do so, we used \emph{quicksort} which is available in MATLAB (R2018a) \citep{MATLAB:2018}. We repeat this process $100$ times. \Cref{fig:plot1} summarizes the results.

\begin{figure}
	\centering
		\setlength{\tabcolsep}{2pt}
	\newcolumntype{C}{>{\centering\arraybackslash} m{0.48\linewidth} }
	\begin{tabular}{m{0.5em}CC}
		&
		$|\dom \underline{P}| = 2^2$
		&
		$|\dom \underline{P}| = 2^6$
		\\
		\rotatebox[origin=l]{90}{$|\mathcal{K}| = 2^6$ \& $|\Omega| = 2^2$}
		&
		\includegraphics[width=\hsize, trim={2cm 0 2cm 0},clip]{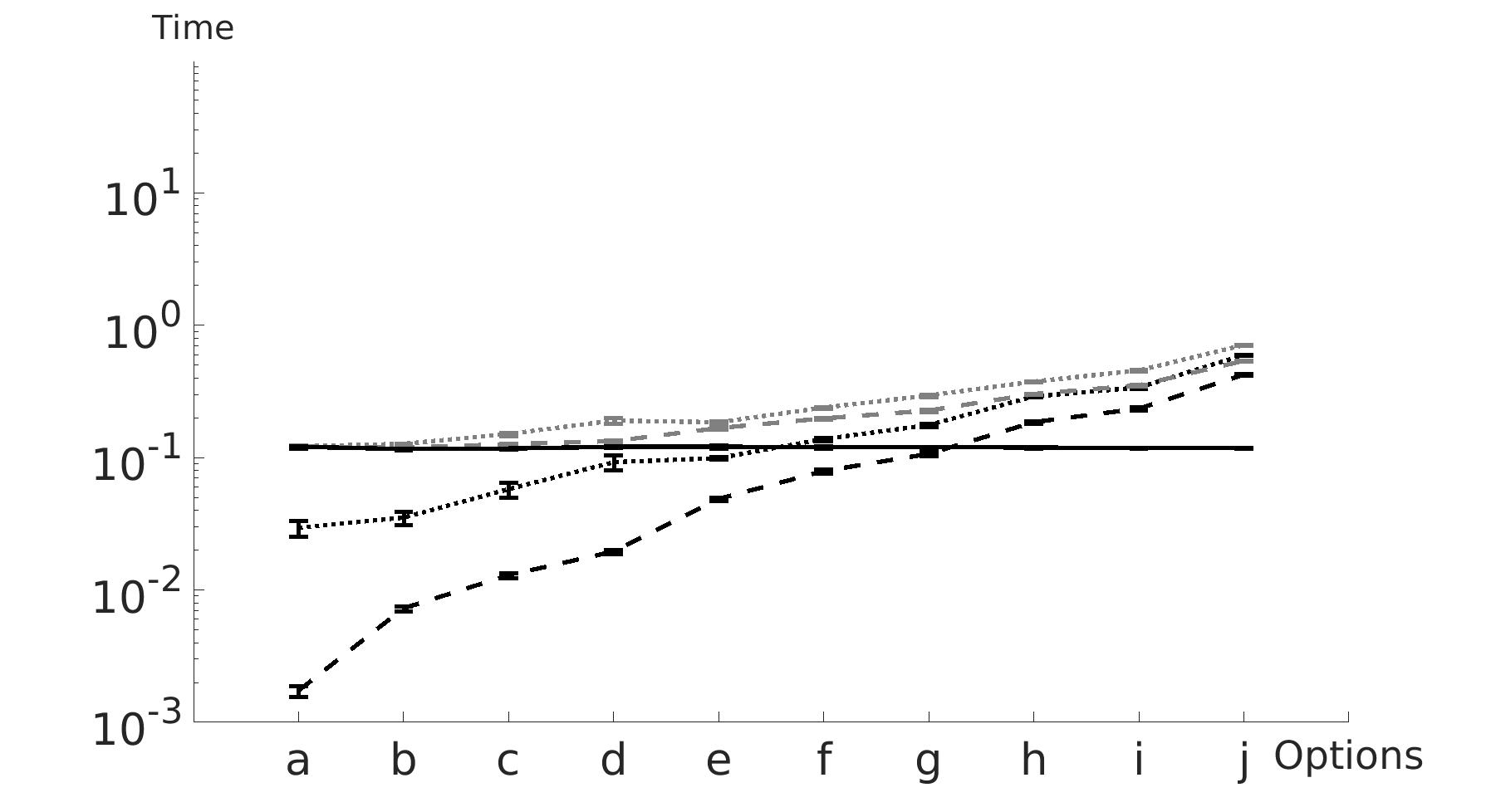}
		&
		\includegraphics[width=\hsize, trim={2cm 0 2cm 0},clip]{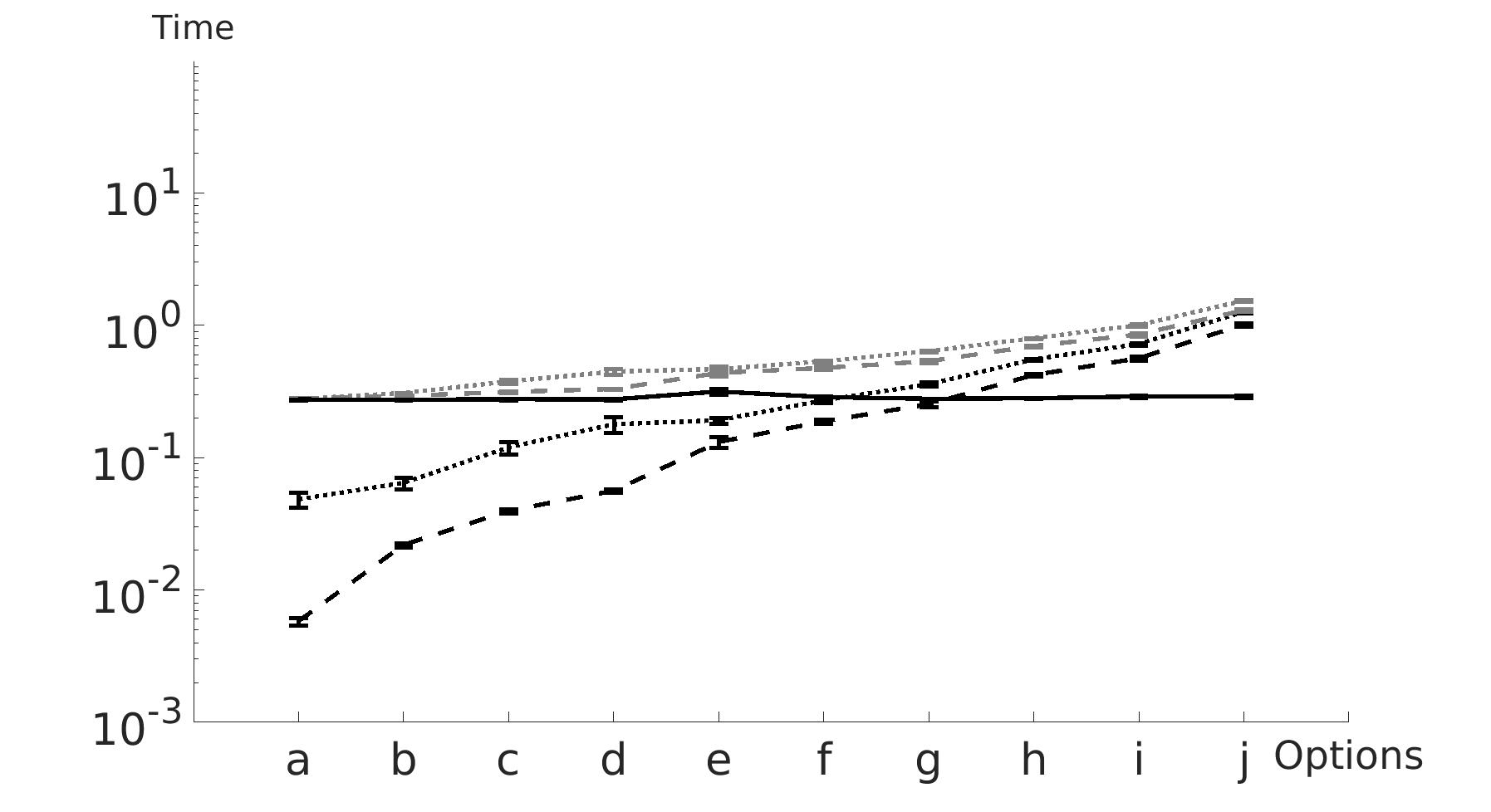}
		\\
		\rotatebox[origin=l]{90}{$|\mathcal{K}| = 2^4$ \& $|\Omega| = 2^6$}
		&
		\includegraphics[width=\hsize, trim={2cm 0 2cm 0},clip]{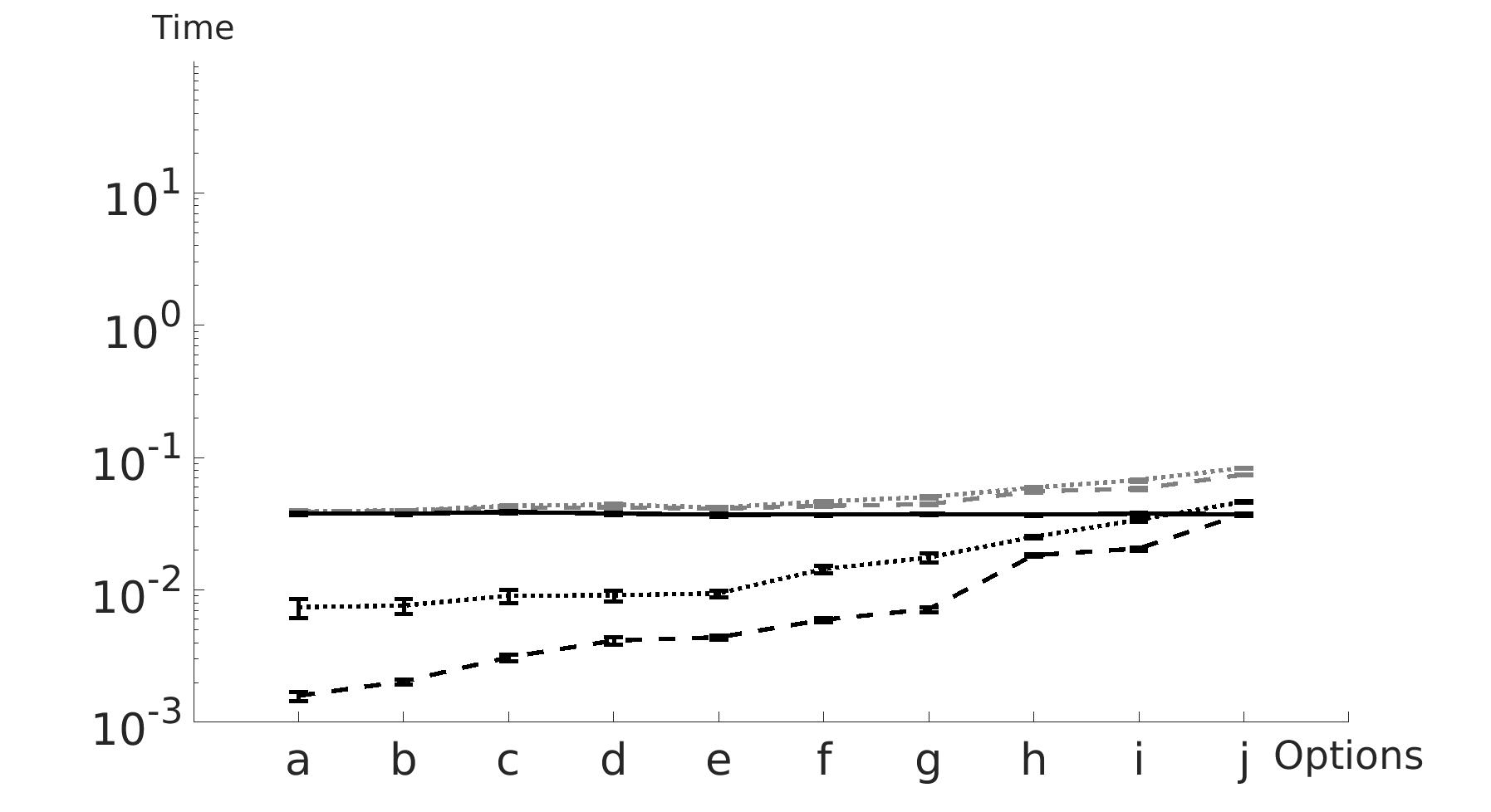}
		&
		\includegraphics[width=\hsize, trim={2cm 0 2cm 0},clip]{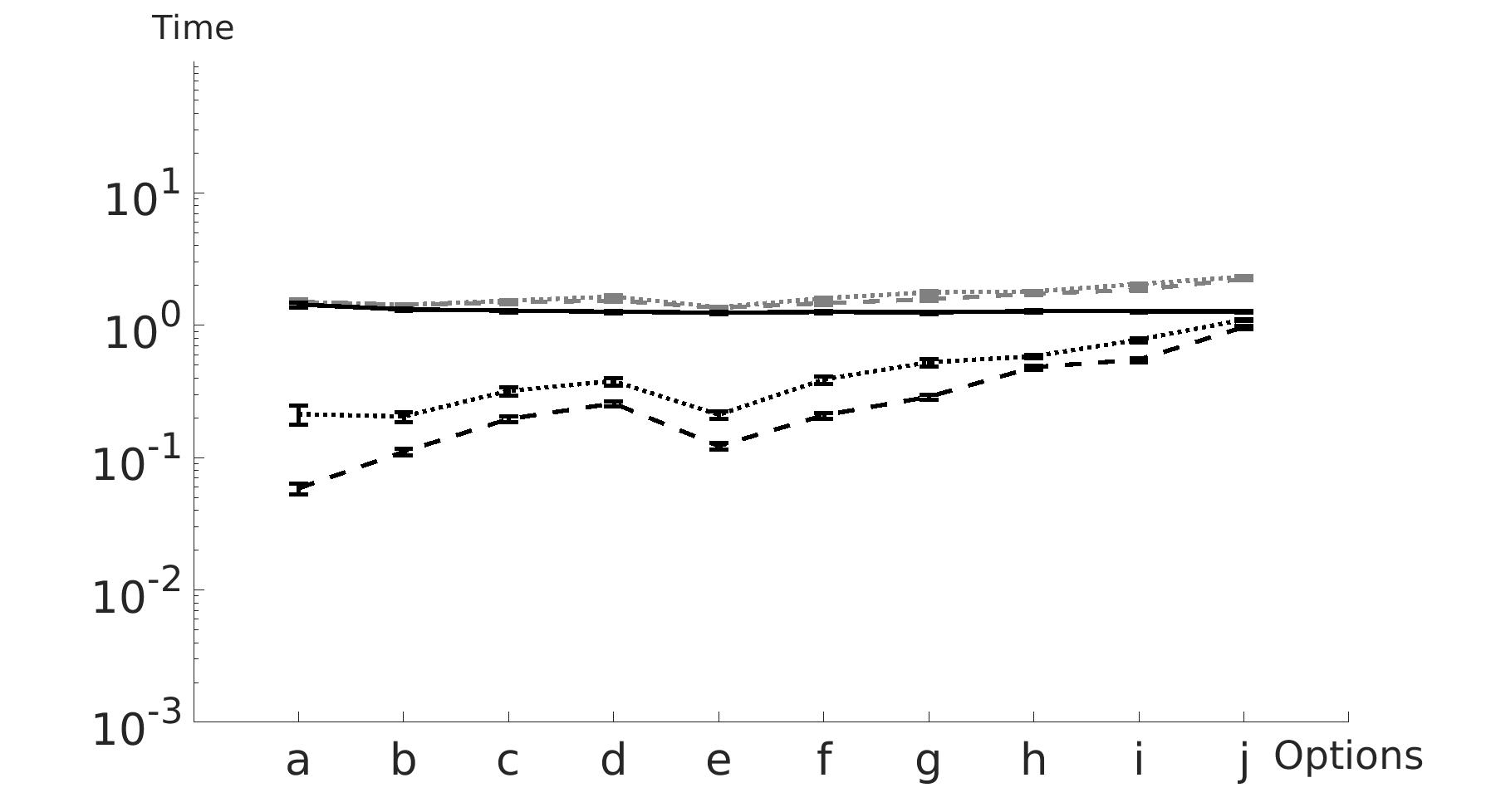}
\end{tabular}
	\includegraphics[width=0.65\linewidth]{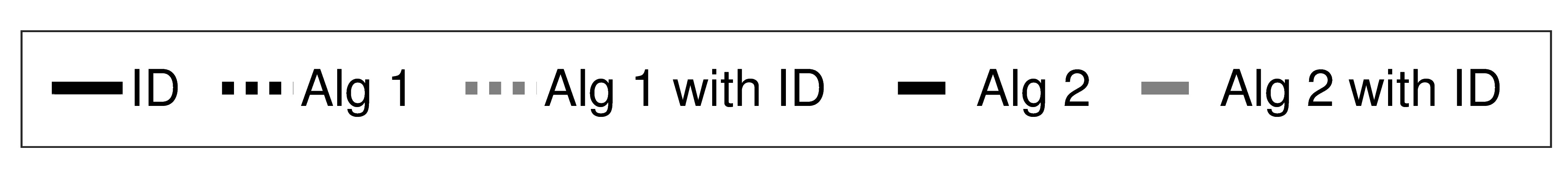}
\caption{Comparison plots of the average computational time for \cref{alg:FindMax:matt,alg:FindMax2} for finding maximal gambles and \cref{alg:FindID} for finding interval dominant gambles. In the left column, $|\dom \underline{P}| = 2^2$ and in the right column, $|\dom \underline{P}| = 2^6$. Each row represents  different numbers of gambles and outcomes with vary options of the numbers of maximal gambles and interval dominant gambles in the set (see \cref{tab:10cases} for each option). The labels indicate algorithms with and without \cref{alg:FindID}.}
\label{fig:plot2}
\end{figure}

\Cref{fig:plot1} shows the average computational time taken during each algorithm with and without \cref{alg:FindID}. The average computational time taken for only \cref{alg:FindID} is also presented there. In the top plots, we show the average computational time for \cref{alg:FindMax:matt,alg:FindMax2,alg:OneLP1,alg:OneLP2}. In the remaining plots, the average computational time for \cref{alg:OneLP1,alg:OneLP2} are so high that their performance are completely dominated by the performance of the other algorithms, and is therefore not presented in these plots. In the left column, the number of outcomes is $2^2$ and in the right column, it is $2^6$. Each row represents a different size of $\mathcal{K}$.

We also consider an impact of the size of $\dom\underline{P}$. \Cref{fig:plot2} shows the average  computational time taken for \cref{alg:FindMax:matt,alg:FindMax2} with and without \cref{alg:FindID}. In the left column, the number of gamble in the domain of $\underline{P}$ is $2^2$ and in the right column, it is $2^6$. Each row represents different numbers of outcomes and gambles.

In both \cref{fig:plot1,fig:plot2}, the horizontal axis indicates different options of $m$, $n$, and $k$ that we consider. The vertical axis presents the computational time which is averaged over 100 random generated sets of gambles. The error bars on the figure represent approximate $95\%$ confidence intervals on the mean computational time.

We also solved \linprogref{P0} by the simplex method in \cref{alg:OneLP1}, using the default simplex method available in MATLAB (R2018a) \citep{MATLAB:2018}. However, this was still slower than \cref{alg:FindMax:matt,alg:FindMax2}. As there is no change in general conclusion, we do not show those plots here.

\section{Discussion and conclusion}\label{sec:conclusion}

In this work, we proposed a new algorithm (\cref{alg:FindMax2}) for finding maximal gambles and compared its performance with \citet[p.~336]{2014:troffaes:itip:computation}'s algorithms (\cref{alg:FindMax:matt}) and \citet{2017:Jansen:Augustin:Schollmeyer} (\cref{alg:OneLP1}). We further improved \cref{alg:OneLP1} by applying the fact that if a gamble is not maximal in one iteration, then it can be excluded from all subsequent iterations (\cref{alg:OneLP2}). We also studied the impact of using interval dominance in \cref{alg:FindID} to eliminate non-maximal gambles as this can reduce the size of the problem.

To find the set of all maximal gambles, \citet{2017:Jansen:Augustin:Schollmeyer}'s algorithm solves a single large linear program for each gamble, while \citet[p.~336]{2014:troffaes:itip:computation}'s algorithm  and our new algorithm solves a larger sequence of smaller linear programs.
For the second case, we proposed early stopping criteria. We also applied common feasible starting points for the entire sequence of problems, based on earlier work \citep{2018:Nakharutai:Troffaes:Caiado}. We found that the primal-dual method can exploit these improvements most effectively, and performs best overall.

To benchmark these algorithms, we presented a new algorithm for generating random sets of gambles with a pre-determined proportion of maximal and interval dominant gambles. This algorithm will be useful for others who want to test their algorithms.
Whilst our benchmarking approach allows careful control over the
properties of the generated decision problems (in particular, the
fraction of optimal gambles according to different decision criteria),
it does have severe computational limitations, due to the need to
evaluate large numbers of natural extensions. Nevertheless, we hope
this work provides a good starting point, and we hope that it inspires
the development of further benchmarking frameworks for testing
algorithms for decision making.

We compared computational performance of  \cref{alg:FindMax:matt,alg:FindMax2,alg:OneLP1,alg:OneLP2} with and without \cref{alg:FindID} on these generated sets.
According to our numerical results, the relative performance of \cref{alg:FindMax:matt,alg:FindMax2,alg:OneLP1,alg:OneLP2} depends on (i) the numbers of outcomes and gambles in the sets, (ii) the ratios of maximal and interval dominant gambles, and (iii) the number of gambles in the domain of lower previsions. If one of these numbers is increasing, then, generally, the average computational time taken on the algorithm is longer. In contrast, the average computational time taken on \cref{alg:FindID} does not depend on the ratios of maximal and interval dominant gambles, but it depends on the numbers of outcomes and gambles in the set and the number of gambles in the domain of lower previsions. This is because  \cref{alg:FindID} has to evaluate the same number of natural extensions regardless of the structure of the problem.

We observed that applying interval dominance (\cref{alg:FindID}) at the beginning benefits \cref{alg:OneLP1,alg:OneLP2} as it makes the linear program smaller, especially if there are many non-interval dominant gambles. Therefore, when using \cref{alg:OneLP1,alg:OneLP2}, we would strongly suggest to run \cref{alg:FindID} first. 
We found that \cref{alg:OneLP2} slightly outperforms \cref{alg:OneLP1}.

In contrast, perhaps surprisingly, interval dominance (\cref{alg:FindID}), at least with our implementation of it, does not help \cref{alg:FindMax:matt,alg:FindMax2}. Even though using \cref{alg:FindID} can eliminate some non-maximal gambles, applying \cref{alg:FindID} first and then performing \cref{alg:FindMax:matt} or \cref{alg:FindMax2} is still slower than performing only \cref{alg:FindMax:matt} or \cref{alg:FindMax2}.
Therefore, we do not recommend applying \cref{alg:FindID} before \cref{alg:FindMax:matt} or \cref{alg:FindMax2}.
This said, there may still be ways to speed up the algorithm for
interval dominance, for example, by adding more stopping criteria.

Overall, both \cref{alg:FindMax:matt,alg:FindMax2} outperform \cref{alg:OneLP1,alg:OneLP2} by an order of magnitude. \Cref{alg:FindMax2} also outperforms \cref{alg:FindMax:matt} in all cases in the experiment, especially when there is only one maximal gamble in the set. In the case that the number of maximal gambles in the set are increasing, there is no big difference in the time taken on \cref{alg:FindMax:matt,alg:FindMax2}, but \cref{alg:FindMax2} still slightly outperforms \cref{alg:FindMax:matt}. When we vary the number of gambles in the domain of lower previsions, the conclusion does not change, i.e., \cref{alg:FindMax2} still outperforms \cref{alg:FindMax:matt}.

Based on theoretical considerations, our newly proposed algorithm, \cref{alg:FindMax2}, is a good choice for implementations as it reduces the number of linear programs, as well as the number of iterations. Our benchmarking study quantified these improvements, and we found that it outperformed all other algorithms tested over all scenarios considered.

\section*{Acknowledgements}

We would like to acknowledge support for this project from Development and Promotion of Science and Technology Talents Project (Royal Government of Thailand scholarship).

\bibliographystyle{plainnat}
\bibliography{references}

\appendix

\section{Proof of technical results  in \cref{sec:benchmark}}\label{pf:appendix}
\begin{proof}[Proof of \Cref{lem:upBoundAlpha}.]
By the definition, we have
\begin{align}
h-\alpha \in \opt_{\succ}(\mathcal{K} \cup \{h - \alpha\})&  \Leftrightarrow \forall g \in \mathcal{K}\colon \overline{E}(h-g -\alpha) \geq 0, \\
& \Leftrightarrow \forall g \in \mathcal{K} \colon \overline{E}(h-g ) \geq \alpha, \\
& \Leftrightarrow \min_{g \in \mathcal{K}} \overline{E}(h - g) \geq \alpha.
\end{align}
Note that 
\begin{equation}
\overline{E}(h - g)  \geq \overline{E}(h - f)- \overline{E}(g-f).
\end{equation}
Suppose that if $g \notin \opt_{\succ}(\mathcal{K})$, then we have $\overline{E}(g-f) < 0$ for some $f \in \opt_{\succ}(\mathcal{K})$, as $g$ is dominated by at least one maximal gamble in $\mathcal{K}$ \citep[p.~336]{2014:troffaes:itip:computation}.
Therefore, for all $g \notin \opt_{\succ}(\mathcal{K})$:
\begin{equation}\label{eq:gIsDomiByf}
\exists f \in  \opt_{\succ}(\mathcal{K}),\ \overline{E}(h - g) \geq \overline{E}(h - f).
\end{equation}
Consequently,
\begin{equation}
h-\alpha \in \opt_{\succ}(\mathcal{K} \cup \{h - \alpha\}) \Leftrightarrow \min_{ f\in  \opt_{\succ}(\mathcal{K})} \overline{E}(h - f) \geq \alpha.
\end{equation}
\end{proof}

\begin{proof}[Proof of  \cref{lem:LowBoundAlpha}.]
	Let $f$ be a maximal gamble in  $\mathcal{K}$. We first show that $f$ is a maximal gamble in  $\mathcal{K}\cup \{h - \alpha\}$ if and only if $\underline{E}(h-f) \leq \alpha$. We see that
	\begin{align}
	f \in \opt_{\succ}(\mathcal{K}\cup \{h-\alpha\})
	& \Leftrightarrow
	\forall g \in \mathcal{K}\cup \{h-\alpha\}\colon
	\overline{E}(f-g) \geq 0 \\
	\intertext{and because $\overline{E}(f-g) \geq 0$ for all $g \in \mathcal{K}$,}
	&\Leftrightarrow
	\overline{E}(f-h+\alpha) \geq 0 \\
	&\Leftrightarrow \overline{E}(f-h) \geq -\alpha\\
	&\Leftrightarrow \underline{E}(h-f) \leq \alpha.
	\end{align}
	Consequently, all maximal gambles in $\mathcal{K}$ are still maximal in $\mathcal{K}\cup \{h - \alpha\}$ if and only if
	\begin{equation}
	\max_{f \in \opt_{\succ}(\mathcal{K})}\underline{E}(h-f) \leq \alpha.
	\end{equation}
\end{proof}

\begin{proof}[Proof of \Cref{lem:max-min_h-f}.]
We first show that \begin{equation}\label{eq:thm1:fg}
\forall f, g \in \opt_{\succ}(\mathcal{K})\colon \underline{E}(h-f) \leq \overline{E}(h-g).
\end{equation} holds.
Let $\mathcal{K}$ be a set of gambles and let $\opt_{\succ}(\mathcal{K})$ be the set of maximal gambles. Suppose that $h$ is another gamble.
Then, for any $f,g \in \opt_{\succ}(\mathcal{K})$, we have
\begin{align}\label{eq:pf:thm1:1}
0 & \leq \overline{E}(f-g) \\ \label{eq:pf:thm1:2}
& = \overline{E}(f-h+h-g) \\ \label{eq:pf:thm1:3}
& \leq \overline{E}(f-h) +\overline{E}(h-g) \\ \label{eq:pf:thm1:4}
& = - \underline{E}(h-f)+\overline{E}(h-g). 
\end{align}
Therefore, by \cref{eq:pf:thm1:4},  for any $f,g \in \opt_{\succ}(\mathcal{K})$:
\begin{equation}
\underline{E}(h-f) \leq \overline{E}(h-g).
\end{equation}
Consequently,
\begin{align}
\max_{f \in \opt_{\succ}(\mathcal{K})}\underline{E}(h-f) &= \underline{E}(h-f^*) \quad \text{for some }  f^* \in  \opt_{\succ}(\mathcal{K})\\
& \leq \overline{E}(h-g) \quad \text{for all }  g \in \opt_{\succ}(\mathcal{K})\quad (\text{by \cref{eq:thm1:fg}})\\
& \leq \min_{ g\in  \opt_{\succ}(\mathcal{K})} \overline{E}(h - g).
\end{align}
Next, we show that 
\begin{equation}
\min_{ f\in  \opt_{\succ}(\mathcal{K})} \overline{E}(h - f) \leq \overline{E}(h) - \max_{g \in \mathcal{K}}\underline{E}(g).
\end{equation}
We see that 
\begin{align}
\min_{ f\in  \opt_{\succ}(\mathcal{K})} \overline{E}(h - f)& = \min_{ g\in  \mathcal{K}} \overline{E}(h - g) \quad (\text{by \cref{eq:gIsDomiByf}})\\
\label{eq:thelastone}  
& \leq \min_{ g\in  \mathcal{K}} \left(\overline{E}(h)-\underline{E}(g) \right)\\
& = \overline{E}(h)- \max_{g \in \mathcal{K}}\underline{E}(g).
\end{align}
\end{proof}

\end{document}